\documentclass{amsart}

\usepackage{amsmath,amsfonts,amssymb,amsthm, multicol}
\usepackage{stmaryrd} 
\usepackage{color} 
\usepackage[all]{xy} 
\usepackage[backref=page]{hyperref} 
\usepackage{enumitem} 

 \AtBeginDocument{%
   \def\MR#1{}
 }

\def\sideremark#1{\ifvmode\leavevmode\fi\vadjust{\vbox to0pt{\vss
\hbox to 0pt{\hskip\hsize\hskip1em%
\vbox{\hsize2cm\tiny\raggedright\pretolerance10000%
\noindent {\color{red}{#1}}\hfill}\hss}\vbox to8pt{\vfil}\vss}}}%

\theoremstyle{plain}
\newtheorem{thm}{Theorem}[section]
\newtheorem{lem}[thm]{Lemma}
\newtheorem{prop}[thm]{Proposition}
\newtheorem{cor}[thm]{Corollary}
\newtheorem*{thm*}{Theorem}

\theoremstyle{definition}
\newtheorem{defi}[thm]{Definition}
\newtheorem{ex}[thm]{Example}

\theoremstyle{remark}
\newtheorem{remi}[thm]{Remark}

\theoremstyle{plain}
\newtheorem{mthm}{Theorem}

\newcommand{\p}{\partial} 
\newcommand{\mL}{\mathcal{L}}  
\newcommand{\bmL}{\boldsymbol{\mathcal{L}}} 

\newcommand{\bZ}{\mathbb{Z}}

\newcommand{\bR}{\mathbb{R}}

\newcommand{\LL}{\left\llbracket} 
\newcommand{\RR}{\right\rrbracket}

\newcommand{\CM}{C^{\infty}(M)}  
\newcommand{\CMG}{C^{\infty}(M)^G}  
\newcommand{\CNH}{C^{\infty}(N)^H}  
\newcommand{\CG}{C^{\infty}(G)}	 
\newcommand{\CmG}{C_m^{\infty}(G)}  
\newcommand{\GA}{\Gamma(A)}  
\newcommand{\GB}{\Gamma(B)}  
\newcommand{\XM}{\mathfrak{X}(M)}  
\newcommand{\XG}{\mathfrak{X}(G)}  
\newcommand{\XmG}{\mathfrak{X}_m(G)}  
\newcommand{\XmH}{\mathfrak{X}_m(H)}  

\newcommand{\CdG}{C_d^{\bullet}(G)}  
\newcommand{\CdH}{C_d^{\bullet}(H)}  
\newcommand{\HdG}{H_d^{\bullet}(G)}  

\newcommand{\CdefG}{C_{\mathrm{def}}^{\bullet}(G)}  
\newcommand{\HdefG}{H_{\mathrm{def}}^{\bullet}(G)}  

\newcommand{\CA}{C^{\bullet}(A)}  
\newcommand{\HA}{H^{\bullet}(A)}  

\newcommand{\CdefA}{C_{\mathrm{def}}^{\bullet}(A)}  
\newcommand{\HdefA}{H_{\mathrm{def}}^{\bullet}(A)}  

\newcommand{\CmbG}{C_m^{\bullet}(G)}  
\newcommand{\CmbH}{C_m^{\bullet}(H)}  
\newcommand{\HCmbG}{H^{\bullet}(C_m^{\bullet}(G))}  
\newcommand{\HCmbH}{H^{\bullet}(C_m^{\bullet}(H))}  

\newcommand{\XmbG}{\mathfrak{X}_m^{\bullet}(G)}  
\newcommand{\XmbH}{\mathfrak{X}_m^{\bullet}(H)}  
\newcommand{\HXmbG}{H^{\bullet}(\mathfrak{X}_m^{\bullet}(G))}  
\newcommand{\HXmbH}{H^{\bullet}(\mathfrak{X}_m^{\bullet}(H))}  

\newcommand{\CmbA}{C_m^{\bullet}(A)}  
\newcommand{\HCmbA}{H^{\bullet}(C_m^{\bullet}(A))}  

\newcommand{\XmbA}{\mathfrak{X}_m^{\bullet}(A)}  
\newcommand{\HXmbA}{H^{\bullet}(\mathfrak{X}_m^{\bullet}(A))}  

\newcommand{\VE}{\mathrm{VE}}  
\newcommand{\oVE}{\overline{\mathrm{VE}}}  

\newcommand{\obullet}{\overline{\bullet}}  
\newcommand{\ocdot}{\overline{\cdot}}  

\newcommand{\Man}{\boldsymbol{\mathrm{Man}}}  
\newcommand{\gpd}{\boldsymbol{\mathrm{gpd}}}  
\newcommand{\St}{\boldsymbol{\mathrm{St}}}  
\newcommand{\Gpd}{\boldsymbol{\mathrm{Gpd}}}  
\newcommand{\GPD}{\boldsymbol{\mathrm{GPD}}}  
\newcommand{\dgla}{\boldsymbol{\mathrm{dgla}}}  
\newcommand{\dglmod}{\boldsymbol{\mathrm{dglmod}}}  
\newcommand{\gla}{\boldsymbol{\mathrm{gla}}}  
\newcommand{\glmod}{\boldsymbol{\mathrm{glmod}}}  

\newcommand{\XmbX}{\mathfrak{X}_m^{\bullet}(\mathcal{X})}
\newcommand{\HXmbX}{H^{\bullet}(\mathfrak{X}_m^{\bullet}(\mathcal{X}))}
\newcommand{\HCmbX}{H^{\bullet}(C_m^{\bullet}(\mathcal{X}))}
\newcommand{\CmbX}{C_m^{\bullet}(\mathcal{X})}
\newcommand{\CX}{C^{\infty}(\mathcal{X})}

\begin{document}

\title{Vector fields and derivations on differentiable stacks}
\date{\today}

\author{Juan Sebastian Herrera-Carmona}
\address{Departamento de Matem\'atica, Universidade Federal do Paran\'a, Setor de Ci\^encias Exatas - Centro Politécnico, 81531-980, Curitiba -Brasil}
\email{sebastianherrera@ufpr.br}

\author{Cristian Ortiz}
\address{Instituto de Matem\'atica e Estat\'istica, Universidade de S\~ao Paulo, Rua do Mat\~ao 1010, Cidade Universit\'aria, 05508-090, S\~ao Paulo - Brasil}
\email{cortiz@ime.usp.br}

\author{James Waldron}
\address{School of Mathematics, Statistics and Physics, Newcastle University, Newcastle upon Tyne, NE1 7RU, UK}
\email{james.waldron@newcastle.ac.uk}

\maketitle

\begin{abstract}
We introduce and study module structures on both the dgla of multiplicative vector fields and the graded algebra of functions on Lie groupoids. We show that there is an associated structure of a graded Lie-Rinehart algebra on the vector fields of a differentiable stack over its smooth functions that is Morita invariant in an appropriate sense. Furthermore, we show that associated Van-Est type maps are compatible with those module structures. We also present several examples.
\end{abstract}


\tableofcontents

\section{Introduction}
The most basic objects in differential geometry are smooth functions and vector fields. 
Associated to a smooth manifold $M$ are the commutative algebra of smooth functions $\CM$ and the Lie algebra of smooth vector fields $\XM$.
The Lie derivative defines a Lie algebra and $\CM$-module isomorphism 

\begin{equation}
\label{eqn: lieder}
\mL : \XM \to \mathrm{Der} (\CM) \; , \; X \mapsto \mL_X
\end{equation}
from $\XM$ to the Lie algebra $\mathrm{Der}(\CM)$ of derivations of $\CM$, and the Lie bracket on $\XM$ satisfies the Leibniz rule
\begin{equation}
\label{eqn: leibniz}
\left[X,fY\right] = \left(\mL_Xf\right) Y + f\left[X,Y\right].
\end{equation}

\noindent In this paper we are concerned with a generalisation of these facts from smooth manifolds to \emph{differentiable stacks}. Recall that a differentiable stack is a stack $\mathcal{X}$ over manifolds admitting a representable epimorphic submersion $M\to \mathcal{X}$ from a manifold to $\mathcal{X}$. As a consequence, the fibered product $M\times_\mathcal{X}M\rightrightarrows M$ with the canonical projections has the structure of a Lie groupoid. Different choices of atlas give rise to Morita equivalent Lie groupoids, establishing a one-to-one correspondence between equivalence classes of differentiable stacks and Morita equivalence classes of Lie groupoids. An immediate observation is that both smooth manifolds and orbifolds can be seen as particular instances of differentiable stacks.

Our main result can be summarised as follows; the explicit formulas and further details will be presented throughout this section.

\begin{thm*}
To each differentiable stack $\mathcal{X}$ there is an associated differential graded Lie algebra $\XmbX$ of vector fields on $\mathcal{X}$ and a differential graded $\XmbX$-module $\CmbX$ of functions on $\mathcal{X}$.
The zeroth cohomology $\CX := H^{0}(\CmbX)$ is a commutative algebra and $\XmbX$ is a complex of $\CX$-modules.\\
\\
The cohomology $\HXmbX$ is a graded $\CX$-module, and the $\XmbX$-module structure on $\CmbX$ descends to cohomology to define a morphism of graded Lie algebras and of graded $\CX$-modules
\[
\bmL : \HXmbX \to \mathrm{Der} (\CX) \; , \; x \mapsto \bmL_x
\]
which is the analogue of \eqref{eqn: lieder}.
The graded Lie bracket $\LL,\RR$ on $\HXmbX$ satisfies the graded Leibniz identity
\[
\LL x,fy \RR = \left(\bmL_{x} f\right) y + (-1)^{|f||x|}f \LL x,y \RR 
\]
for appropriate $x,y,f$, which is the analogue of \eqref{eqn: leibniz}.
Up to an appropriate notion of quasi-isomorphism these objects are all independent of the choices involved in the construction and are unchanged if $\mathcal{X}$ is replaced by an equivalent stack $\mathcal{Y}$.\\
\\
If $\mathcal{X}$ is isomorphic to a manifold $M$ then $\XmbX$ is quasi-isomorphic to $\XM$, $\CmbX$ is quasi-isomorphic to $\CM$, and the various module structures are equivalent to the standard ones.
\end{thm*}

\begin{remi}
Since the algebra of functions $\CX$ is concentrated in degree zero, the sign appearing in the previous graded Leibniz identity is actually positive. See Definition \ref{def: grLR}.
\end{remi}

\begin{remi}
\label{rem: non-prop}
If $\mathcal{X}$ is proper, i.e.\ if $\mathcal{X}$ is equal to the quotient stack $M//G$ associated to a proper Lie groupoid $G\rightrightarrows M$, then the cohomology $\HCmbX$ is equal to $\CX$ placed in degree zero.
In the non-proper setting $\HCmbX$ is a graded algebra and $\HXmbX$ is a graded $\HCmbX$-module, but the situation is more complicated, see \S \ref{ex: non-prop}.
\end{remi}

\begin{remi}
The differential graded Lie algebra $\XmbX$ was constructed independently in work of Berwick-Evans \& Lerman \cite{BerwickEvansL20} and of the second two named authors \cite{OrtizW19}, and the underlying cochain complex of $\CmbX$ is a truncation of the differential groupoid complex introduced in \cite{Crainic03}.
Our new contributions in this work are: establishing the $\XmbX$-module structure on $\CmbX$ and its properties; explaining the relation of $\CmbX$ to other notions of smooth functions on differential stacks; the construction of a certain infinitesimal model of these objects and structures, valid whenever certain connectivity assumptions on a groupoid presentation are satisfied. We also present several examples throughout the paper.
\end{remi}

We can generalise the right hand side of \eqref{eqn: lieder}.
There are a number of equivalent definitions of the algebra $C^{\infty}(\mathcal{X})$ of functions on a differentiable stack $\mathcal{X}$, some of which depend on the choice of an atlas $M \to \mathcal{X}$ with associated Lie groupoid $G$ over $M$.
These include:
\begin{enumerate}
\item[(i)] stack morphisms from $\mathcal{X}$ to the manifold $\bR$, 
\item[(ii)] global sections of the structure sheaf $\mathcal{O}_{\mathcal{X}}$ of $\mathcal{X}$, 
\item[(iii)] $G$-invariant functions on $M$, 
\item[(iv)] morphisms of Lie groupoids from $G$ to the unit Lie groupoid $\bR \rightrightarrows \bR$.
\end{enumerate}

In order to further capture the ``stacky'' nature of $\mathcal{X}$, our first main result compares certain extensions and variations of these definitions, each of which involves moving from sets to \emph{categories} or from vector spaces to \emph{cochain complexes}:
\begin{enumerate}
\item[(i)]  The \emph{category} $\mathrm{Hom}_{\St}(\mathcal{X},B\bR)$ of morphisms from $\mathcal{X}$ to the classifying \emph{stack} $B\bR$ of the Lie group $\bR$; 
\item[(ii)]  The total sheaf cohomology $H^{\bullet}(\mathcal{O}_{\mathcal{X}})$, which carries a natural graded algebra structure \cite[\S 3]{BehrendX03}; 
\item[(iii)]  The differentiable groupoid cochain complex $\CdG$, which carries a natural differential graded algebra structure \cite[\S 1.2]{Crainic03}; 
\item[(iv)]  The \emph{category} $\mathrm{Hom}_{\Gpd}(G,\bR)$ of morphisms from $G$ to the Lie groupoid $\bR \rightrightarrows \ast$ where $\bR$ is considered as an abelian Lie group.
\end{enumerate}
There are some known relations between these objects: if $M$ is Hausdorff and paracompact then the cohomology $\HdG$ of $\CdG$ is isomorphic to $H^{\bullet}(\mathcal{O}_{\mathcal{X}})$ as a graded algebra \cite[\S 3]{BehrendX03}, and the 1-cocycles $Z^1_d(G) \subseteq \CdG$, which we will denote by $\CmG$, coincide with the morphisms from $G$ to the Lie group $\bR$, which are also called \emph{multiplicative functions} on $G$ \cite{MackenzieX98,MackenzieBook05}.

To formulate the following Theorem we need two further objects. 
We denote by $\CmbG$ the 2-term truncation $C^{\infty}(M) \xrightarrow{\delta} \CmG$ of the cochain complex $\CdG$.
Associated to $\CmbG$ is a category $\CM \ltimes \CmG$ in which the objects are multiplicative functions $F \in \CmG$ and a morphism $F \to F'$ is identified with a function $f \in \CM$ such that $F' = F + \delta f$, see \S \ref{sec: 2vec}.

\begin{mthm}
[= Theorem \ref{thm: fun}]
\label{mthm: fun}
Let $G\rightrightarrows M$ be a Lie groupoid.
The following two categories are isomorphic:
\begin{enumerate}
\item 
The category $\CM \ltimes \CmG$ associated to the 2-term complex $\CmbG$.
\item 
The category $\mathrm{Hom}_{\Gpd} (G,\bR)$ of Lie groupoid morphisms from $G$ to $\bR$ considered as a Lie groupoid over a point.
\end{enumerate}
Additionally, if $M$ is Hausdorff and paracompact then \eqref{it: thm_fun_CmbG} and \eqref{it: thm_fun_HomGpd} are equivalent to:
\begin{enumerate}
\item[(c)] 
The category $\mathrm{Hom}_{\St} (BG,B\bR)$ of morphisms of stacks from $BG$ to $B\bR$.
\end{enumerate}
More generally, if $\mathcal{X}$ is a differentible stack, $M \to \mathcal{X}$ is an atlas with $M$ Hausdorff and paracompact, and $G$ is the associated Lie groupoid over $M$, then the categories \eqref{it: thm_fun_CmbG} and \eqref{it: thm_fun_HomGpd} are equivalent to:
\begin{enumerate}
\item[(d)] The category $\mathrm{Hom}_{\St} (\mathcal{X},B\bR)$ of morphisms of stacks from $\mathcal{X}$ to $B\bR$.
\end{enumerate}
\end{mthm}

\begin{remi}
Theorem \ref{mthm: fun} can be applied to \emph{any} differentiable stack: if $M \to \mathcal{X}$ is an arbitrary atlas then by replacing $M$ with the disjoint union of a suitable open cover if necessary, one can arrange that $M$ is Hausdorff and paracompact.
\end{remi}

Our second main result involves an action of a certain differential graded Lie algebra $\XmbG$ on the complex $\CmbG$.
In \cite{Hepworth09} Hepworth defined the category $\Gamma(T\mathcal{X})$ of sections of the tangent stack $T\mathcal{X}$ of a differentiable stack $\mathcal{X}$ and showed that this category is equivalent to the category $\Gamma(TG)$ of multiplicative sections of the tangent groupoid of a Lie groupoid $G$ presenting $\mathcal{X}$.
It was shown independently by Berwick-Evans \& Lerman \cite{BerwickEvansL20} and by the second two authors \cite{OrtizW19} that the category $\Gamma(TG)$ carries a natural structure of strict Lie 2-algebra, which, up to a natural notion of equivalence, is independent of the choice of the Lie groupoid $G$.

In terms of differential graded Lie algebras, or dgla's for short, the results in \cite{BerwickEvansL20} and \cite{OrtizW19} show that associated to a Lie groupoid $G$ there is a dgla $\XmbG$ concentrated in degrees -1 and 0, and if $H$ is a Lie groupoid Morita equivalent to $G$ then the dgla's $\XmbG$ and $\XmbH$ are quasi-isomorphic.
The dgla $\XmbG$ is constructed using the Lie algebroid of $G$ and the Lie algebra of multiplicative vector fields on $G$, which are vector fields $X : G \to TG$ that are Lie groupoid morphisms. 
See \ref{sec: XmbG} for the precise definitions.
We note that \cite{BonechiCLGX20} has generalised some of this to multi-vector fields.

\begin{mthm}
[= Theorems \ref{thm: mu} and \ref{thm: morita}]
\label{mthm: mu}
The following statements hold.
\begin{enumerate}
\item If $G$ is a Lie groupoid then the map 
\begin{align*}
\mu: \XmbG \otimes \CmbG & \to \CmbG \\
(\alpha , X) \otimes (f , F) & \mapsto (u^*\mL_{\alpha^r} F + \mL_{X_M} f , \mL_{X} F)
\end{align*}
defines a left dg-Lie module structure of $\XmbG$ on $\CmbG$.

\item \label{it: mthm_morita_GH}
If $G$ and $H$ are Morita equivalent Lie groupoids then $(\XmbG,\CmbG)$ and $(\XmbH , C_m^{\bullet}(H))$ are quasi-isomorphic objects in the category $\boldsymbol{\mathrm{dglmod}}$.

\item \label{it: mthm_morita_X}
If $\mathcal{X}$ is a differentiable stack and $M \to \mathcal{X}$ is an atlas then there is an associated object $(\XmbG,\CmbG)$ in $\boldsymbol{\mathrm{dglmod}}$ where $G$ is the Lie groupoid associated to $M$.
Up to quasi-isomorphism in $\boldsymbol{\mathrm{dglmod}}$ this object does not depend on the choice of atlas of $\mathcal{X}$.
\end{enumerate}
\end{mthm}

In Theorem \ref{mthm: mu}, $X$ is a multiplicative vector with associated vector field $X_M$ on $M$, $\alpha^r$ is the right-invariant vector field associated to the section $\alpha$ of the Lie algebroid of $G$, $f \in \CM$, $F \in \CmG$ is a multiplicative function, and $u: M \to G$ is the unit map of $G$.
See \S \ref{sec: mod} for the precise definitions.
We note that the operator $u^*\mL_{\alpha^r}$ also appears in \cite[\S 4]{AbadC11}.
The category $\dglmod$ and the notion of quasi-isomorphism therein is defined in Definition \ref{def: dglmod}.

As a consequence, one has the following. Given a differentiable stack $\mathcal{X}$ and $G$ a Lie groupoid whose classifying stack is isomorphic to $\mathcal{X}$, then:

\begin{mthm}
[= Theorem \ref{thm: co_stacks}]
There is an associated a graded Lie-Rinehart algebra 
\[
\mathfrak{X}^{\bullet}(\mathcal{X}) := \HXmbG
\]
over the algebra $\CX$, whose isomorphism class is independent of the choice of atlas. Also, if $\mathcal{X}$ is equivalent to $\mathcal{Y}$, then $\mathfrak{X}^{\bullet}(\mathcal{X})\cong \mathfrak{X}^{\bullet}(\mathcal{Y})$ as graded Lie-Rinehart algebras over $\CX$.
\end{mthm}

One can bring Lie theory into the picture by considering a Lie groupoid $G$ with Lie algebroid $A$. In this case, there are a complex $\mathfrak{X}^{\bullet}_{m}(A)$ of infinitesimally multiplicative functions on $A$ (Definition \ref{def: CmbA}), as well as a dgla $\mathfrak{X}^{\bullet}_{m}(A)$ of derivations of $A$ (Definition \ref{def: XmbA}), together with certain module structures introduced in Definitions \ref{def: obull} and \ref{def: ocdot}, respectively. We show that there are Van-Est maps compatible with the module structures as in the next result. 

\begin{mthm}
[= Theorem \ref{thm: VE}]
\label{mthm: VE}
The following statements hold:
\begin{enumerate}
\item 
The Van-Est map $\VE : \XmbG \to \XmbA$ is a morphism of dgla's.
\item The following diagram is a commutative diagram of morphisms of cochain complexes:
\begin{equation}
\label{diag: mVE_mod}
\xymatrix{
\XmbG \otimes \CmbG \ar[d]_{\VE \otimes \oVE} \ar[r]^-{\mu} & \CmbG \ar[d]^{\oVE} \\
\XmbA \otimes \CmbA \ar[r]_-{\overline{\mu}} & \CmbA
}
\end{equation}
\item The vertical arrows \eqref{diag: mVE_mod} are isomorphisms whenever $G$ is source simply connected.
\end{enumerate}
\end{mthm}

\begin{remi}
Instead of the (differential graded) algebras $\CmbG$ and $C^{\infty}(\mathcal{X})$, one can study the in-general noncommutative convolution algebra $\mathcal{A}_{G}$ of $G$.
It is shown in \cite{BjarneP2023} that multiplicative vector fields act as derivations on $\mathcal{A}_{G}$, and this action extends to a morphism of cochain complexes
\[
C_{\mathrm{def}}^{\ge 1}(G) \to C^{\bullet}(\mathcal{A}_{G},\mathcal{A}_{G})
\]
from a truncation of deformation complex of $G$ to the Hochschild cohomology of $\mathcal{A}_{G}$, see loc.\ cit.\ for further details.
\end{remi}

The paper is organized as follows. In \S 2 we present the necessary background in order to state and prove our main results. In \S 3 we recall the complex of multiplicative functions on a Lie groupoid as well as the dgla of multiplicative vector fields. We also introduce the infinitesimal version of these objects and we show the existence of a Van-Est map relating them. \S 4 is devoted to the study of several notions of functions on differentiable stacks and their relation to the complex of multiplicative functions on a Lie groupoid presenting a given a stack. In \S 5 and \S 6 we introduce and study the module structures appearing in the main results explained above. In \S 7 we show the graded Lie-Rinehart algebra structure on vector fields on a differentiable stack. \S 8 explains the compatibility of the Van-Est maps with the module structures introduced in \S 5 and \S 6. In \S 9 we present several examples.

\subsection{Acknowledgements}
C. Ortiz thanks the Newcastle University for the hospitality while part of this work was being done. The research of C. Ortiz was partially supported by the National Council of Research and Development CNPq-Brazil, Bolsa de Produtividade em Pesquisa Grant 315502/2020-7 and by Grant 2016/01630-6 Sao Paulo Research Foundation - FAPESP.

\section{Background and notation}
All vector spaces, cochain complexes, algebras and other algebraic structures will be over the real numbers $\bR$.

By a manifold we will always mean a finite dimensional smooth manifold Hausdorff and paracompact.
The symbols $M$ and $N$ will always denote manifolds.
We denote by $\CM$ the algebra of smooth functions on $M$ and by $\XM$ the Lie algebra of vector fields on $M$.
The derivative of a smooth map $\phi : M \to N$ is denoted by $d\phi$.
The symbols $X$ and $Y$ will always denote vector fields.
If $\phi : M \to N$ is a diffeomorphism and $X \in \XM$ then we denote by $\phi_*X$ the vector field on $N$ defined by $y \mapsto d\phi(X(\phi^{-1}(y)))$.
We note that if $f \in C^{\infty}(N)$ then $\phi^* (\mL_{\phi_* X} f) = \mL_X (\phi^* f)$.

\subsection{Lie groupoids}
\label{sec: lie_groupoids}
For general facts about Lie groupoids see \cite{MackenzieBook05}, and \cite[\S 5.4]{MoerdijkMBook03} or \cite[\S 2.1\&2.5]{MoerdijkMBook05} for the notion of weak and Morita equivalence.
Unless otherwise stated, the symbol $G$ will always denote a Lie groupoid with base $M$.
We denote the structure maps of $G$ by $s,t : G \to M$ (source and target), $u: M \to G$ (unit), $m : G_2 \to G$ (multiplication) and $i: G \to G$ (inversion), where $G_2 := G \times_{M} G$ is the manifold of composable pairs of morphisms in $G$.
If $g \in G$ then the associated right and left translations are the diffeomorphisms $R_g : s^{-1}(t(g)) \to s^{-1}(s(g))$, $h \mapsto hg$, and $L_g : t^{-1}(s(g)) \to t^{-1}(t(g))$, $h \mapsto gh$, respectively.

We denote the strict 2-category of Lie groupoids by \textbf{Gpd}.
If $G$ and $H$ are Lie groupoids then the category $\mathrm{Hom}_{\Gpd}(G,H)$ has objects the Lie groupoid morphisms $\phi : G \to H$, and arrows $\psi : \phi \Rightarrow \phi'$ given by smooth natural transformations.

\subsection{Lie algebroids}
\label{sec: lie_algebroids}
For the basics on Lie algebroids see \cite{MackenzieBook05} or \cite{MoerdijkMBook03}. 
We define the Lie algebroid $A$ of the Lie groupoid $G$ by $A := (\mathrm{Ker} \, ds) |_{M}$ where the restriction to $M$ is via pullback along the unit map $u: M \to G$.
(Note that some references define $A$ to be equal to $(\mathrm{Ker} \, dt) |_{M}$.)
The symbols $\alpha$ and $\beta$ will always denote sections of $A$.

The right-invariant and the left-invariant vector fields on $G$ associated to $\alpha \in \Gamma(A)$ are defined by
$\alpha^{r}(g) = dR_g (\alpha(1_{t(g)}))$ and $\alpha^l(g) = - dL_g (di (\alpha(1_{s(g)})))$, respectively.
Note that in some references the definition of $\alpha^l$ differs from ours by a sign.
With our convention $\alpha^l = - i_* \alpha^r$.
The map $\alpha \mapsto \alpha^r$ (resp.\ $\alpha \mapsto \alpha^l$) is a vector space isomorphism from $\GA$ to the Lie algebra of right (resp.\ left) invariant vector fields on $G$ and the Lie bracket on $\Gamma(A)$ is defined such that the first of these maps is a Lie algebra isomorphism.
The anchor map $a : A \to TM$ is the restriction to $M$ of the map $dt : \mathrm{Ker} \, ds \to t^* TM$.

\subsection{Differentiable stacks}
\label{sec: stacks}
The main properties of differentiable stacks can be found in \cite{BehrendX03}.
By a stack we always mean a pseudo-functor $\Man \to \gpd$ from the category $\Man$ of smooth manifolds to the 2-category $\gpd$ of set-theoretic groupoids.
We denote the bicategory of stacks by $\St$.
A \emph{differentiable stack} $\mathcal{X}$ is a stack for which there exists a representable epimorphic submersion $M \to \mathcal{X}$ from a manifold $M$ to $\mathcal{X}$.
In this case the manifold $G := M \times_{\mathcal{X}} M$ carries a natural Lie groupoid structure over $M$, $\mathcal{X}$ is equivalent to the classifying stack $BG$ of principal $G$-bundles, and we say that the Lie groupoid $G$ presents $\mathcal{X}$.
Different atlases lead to Morita equivalent Lie groupoids, and more generally two Lie groupoids are Morita equivalent if and only if the corresponding classifying stacks are equivalent.

\subsection{Differential graded objects}
\label{sec: dgalgebra}
For background material on differential graded Lie algebras see \cite[\S IV,V]{Manetti04}.
By complex we will always mean cochain complex of real vector spaces, so that all differentials have degree $+1$.
We will use the abbreviations `dga' (differential graded algebra), `dgla' (differential graded Lie algebra), and `gla' (graded Lie algebra).
A gla is a dgla with zero differential.
We denote the category of dgla's by $\dgla$ and the subcategory of gla's by $\gla$.
We denote by $\mathrm{End}^{\bullet} (C^{\bullet})$ the dgla of graded endomorphisms of a cochain complex $C^{\bullet}$, and by $\mathrm{Der}^{\bullet} (B^{\bullet})$ the dgla of graded derivations of a dga $B^{\bullet}$.

A \emph{differential graded Lie module structure} (or just a dgla $L^{\bullet}$-module for short) over a dgla $(L^{\bullet},[\cdot,\cdot],\partial)$ is a complex $(M^{\bullet},\delta)$ equipped with a degree zero morphism of dgla's $L^{\bullet} \to \mathrm{End}^{\bullet}(M^{\bullet})$, or equivalently a degree zero cochain map $L^{\bullet} \otimes M^{\bullet} \to M^{\bullet}$, $x \otimes y \mapsto x \bullet y$ satisfying 
\begin{equation}\label{eqn:dglamod}
\begin{aligned}
  & \delta(x\bullet y)=\partial(x)\bullet y+(-1)^{|x|}x\bullet \delta(y), \quad \text{and} \\
  & [x,x'] \bullet y = x \bullet (x' \bullet y) - (-1)^{|x||x'|} x' \bullet (x \bullet y);
\end{aligned}
\end{equation}
for all homogeneous elements $x, x' \in L^{\bullet}$ and $y \in M^{\bullet}$.
If $\phi : L^{\bullet} \to K^{\bullet}$ is a morphism of dgla's and $M^{\bullet}$ is a $K^{\bullet}$-module then the \emph{pullback module} is the $L^{\bullet}$-module $\phi^* M^{\bullet}$ with underlying cochain complex $M^{\bullet}$ and action $x \bullet y := \phi(x) \bullet y$.
If $L^{\bullet}$ is a gla then a \emph{graded module or dg-module} over $L^{\bullet}$ is a differential graded module $C^{\bullet}$ for which the differential on $C^{\bullet}$ is equal to zero.

\begin{defi}
\label{def: dglmod}
The category $\dglmod$ has objects given by pairs $(L^{\bullet},M^{\bullet})$ where $L^{\bullet}$ is a dgla and $M^{\bullet}$ is a $L^{\bullet}$-module, a morphism $(L^{\bullet},M^{\bullet}) \to (L'^{\bullet},M'^{\bullet})$ is a pair $(\phi,\psi)$, where $\phi: L^{\bullet} \to L'^{\bullet}$ is a morphism of dgla's and $\psi : \phi^* M' \to M$ is a morphism of $L^{\bullet}$-modules, and the composition of morphisms of given componentwise.
\end{defi}

A \emph{quasi-isomorphism} in $\dglmod$ is a morphism $(\phi,\psi)$ in which $\phi$ is a quasi-isomorphism of dgla's and $\psi$ is a quasi-isomorphism of cochain complexes.
A pair of objects $(L^{\bullet},M^{\bullet})$ and $(L'^{\bullet},M'^{\bullet})$ are \emph{quasi-isomorphic} if they isomorphic in the category obtained by inverting all quasi-isomorphisms in $\dglmod$.

\begin{remi}
\label{rem: qisom_identify}
Suppose that $(\phi,\psi) : (L^{\bullet},M^{\bullet}) \to (L'^{\bullet},M'^{\bullet})$ is a quasi-isomorphism in $\dglmod$.
If we identify $H^{\bullet}(L'^{\bullet})$ with $H^{\bullet}(L^{\bullet})$ via the isomorphism $H(\phi)$ then $H^{\bullet}(\psi) : H^{\bullet}(M'^{\bullet}) \to H^{\bullet}(M^{\bullet})$ is an isomorphsim of graded $H^{\bullet}(L^{\bullet})$-modules.
\end{remi}

\begin{defi}
\label{def: glmod}
The category $\boldsymbol{\mathrm{glmod}}$ is defined as follows: the objects are pairs $(L^{\bullet},M^{\bullet})$ where $L^{\bullet}$ is a gla and $M^{\bullet}$ is a graded module over $L^{\bullet}$, a morphism $(L^{\bullet},M^{\bullet}) \to (L'^{\bullet},M'^{\bullet})$ is a pair $(\phi,\psi)$, where $\phi: L^{\bullet} \to L'^{\bullet}$ is a morphism of gla's and $\psi : \phi^* M' \to M$ is a morphism of graded $L^{\bullet}$-modules, and the composition of morphisms of given componentwise.
\end{defi}

It follows from Remark \ref{rem: qisom_identify} that there is a natural functor $\boldsymbol{\mathrm{dglmod}} \to \boldsymbol{\mathrm{glmod}}$ mapping objects to their cohomology.
This functor maps quasi-isomorphisms in $\boldsymbol{\mathrm{dglmod}}$ to isomorphisms in $\boldsymbol{\mathrm{glmod}}$.

\subsection{2-term complexes and 2-vector spaces}
\label{sec: 2vec}
Every cochain complex $C^{\bullet} = C^0 \xrightarrow{\delta} C^1$ concentrated in degrees 0 and 1 determines a groupoid $C^0 \ltimes C^1\rightrightarrows C^1$, where
a pair $(x,y) \in C^0 \times C^1$ is by definition a morphism from $y$ to $y + \delta x$, and the composition is given by $(x',y') \circ (x,y) = (x+x',y)$, defined whenever $y' = y + \delta x$.
The sets of objects and morphisms are each vector spaces, and the structure maps of the category are all linear, so that $C^0 \ltimes C^1$ is a category internal to the category of vector spaces, or a \emph{2-vector space} for short.
This construction extends in an evident way to complexes concentrated in degrees $i$ and $i+1$ for some fixed $i \in \bZ$.
Restricting to the $i=0$ case the assignment
\[
C^0 \xrightarrow{\delta} C^1 \mapsto C^0 \ltimes C^1 \rightrightarrows C^1
\]
extends to an equivalence between the (suitably defined) 2-category of 2-term cochain complexes and the 2-category of 2-vector spaces, see \cite[\S 3]{BaezC04} or \cite[\S 2]{Daenzer14}.
(The original reference is \cite[Expos\'e XVIII, \S 1.4]{SGA4III73} in which 2-vector spaces are called `Picard groupoids'.)
This can also be seen as a special case of the Dold-Kan correspondence \cite[\S 8.4]{WeibelBook94}.
Note that there are other unrelated notions of 2-vector spaces, e.g.\ \cite{KapranovV91}.

\subsection{dgla's, Lie 2-algebras and crossed modules}
\label{sec: 2lie}
The equivalence of 2-term complexes and 2-vector spaces described in \S \ref{sec: 2vec} extends to the setting of Lie algebras.
If $L^{\bullet}$ is a dgla concentrated in degrees -1 and 0 then the 2-vector space $L^{-1} \ltimes L^0 \rightrightarrows L^0$ associated to the underlying complex $L^{-1} \to L^0$ is in a natural way a category internal to the category of Lie algebras or \emph{strict Lie 2-algebra} for short. That is,
$L^0$ is a Lie algebra, the vector space $L^{-1} \times L^0$ carries the semi-direct product Lie bracket defined by the adjoint action of $L^0$ on $L^1$, and the structure maps are all Lie algebra homomorphisms.

As shown in \cite[Def.\ 47]{BaezC04}, 2-term dgla's are also equivalent to \emph{crossed modules of Lie algebras}. Concretely, a dgla $L^{\bullet}$ concentrated in degrees -1 and 0 defines the crossed module
\[
L^{-1} \xrightarrow{\delta} L^0 \xrightarrow{\mathrm{ad}} \mathrm{Der}(L^{-1})
\]
where the Lie algebra maps $\delta$ and $\mathrm{ad}$ are given by the differential of $L^{\bullet}$ and the restriction of the adjoint representation, respectively. See \cite[\S 5]{BaezC04} for further details.

\subsection{Graded Lie-Rinehart algebras}
\label{sec: grLR}
Lie-Rinehart algebras, first introduced in \cite{Rinehart63}, are the algebraic analogues of differential-geometric Lie algebroids. Note that in some references Lie-Rinehart algebras are called Lie algebroids, but we will keep separate the terminology to avoid confusion.
Differential graded Lie-Rinehart algebras appear in a number of references, including \ \cite[\S 2.2]{Grivaux20}, \cite[\S 2.1]{Nuiten19}, \cite[\S 2]{Vezzosi15}, \cite{Vitagliano2014}. In this work we will only need a simpler notion of a graded Lie-Rinehart algebra over an \textbf{$\bR$-algebra} as given in the following definition.

\begin{defi}
\label{def: grLR}
A \emph{graded Lie-Rinehart algebra} over an $\mathbb{R}$-algebra $R$ consists of a graded Lie algebra $L^\bullet$, a graded left $R$-module structure on the underlying graded vector space of $L^{\bullet}$, and a morphism of both gla-modules and graded $R$-modules $a : L^\bullet \to \mathrm{Der} (R)$. 
This data must satisfy the graded Leibniz rule
\begin{equation}
\label{eqn: grLeib}
[x,ry] = \left(a(x) (r)\right) y + r[x,y]
\end{equation}
for every $r \in R$ and homogeneous elements $x,y \in L^{\bullet}$.
\end{defi}

\begin{remi}Note that the Leibniz identity in Definition \ref{def: grLR} does not involve a minus sign of the form $(-1)^{|r||x|}$. Indeed, since the algebra $R$ is concentrated in degree zero, then $\mathrm{Der} (R)$ is also concentrated in degree zero, hence a Lie algebra. 
\end{remi}
\begin{defi}
\label{def: gLRalg} The category $\mathbf{gLRalg}$ has object pairs $(L^{\bullet}, R)$ consisting of a graded Lie-Rinehart algebra $L^{\bullet}$ over $R$. A morphism $(L^{\bullet}, R) \to (L'^{\bullet}, R')$ is a pair $(\phi, \psi)$ where $\psi: R \to R'$ is a morphism of $\mathbb{R}$-algebras, and $\phi: L^{\bullet} \to L'^{\bullet}$ is both a morphism of graded Lie algebras and a morphism of $R$-modules. Here, $L'$ is an $R$-module with the pullback module structure induced by $\psi: R \to R'$.
\end{defi}


\section{Some 2-term complexes and dgla's}
\label{sec:bg_notation}
In this section we define several objects naturally associated to a Lie groupoid $G$ and its Lie algebroid $A$. There are two dgla's: 
\begin{itemize}
\item $\XmbG$: the dgla of \emph{multiplicative vector fields} on $G$,
\item $\XmbA$: the dgla of \emph{Lie algebroid derivations} of $A$.
\end{itemize}
Also there are two cochain complexes, each of them carrying a graded algebra structure on their cohomology:
\begin{itemize}
\item $\CmbG$: the complex of \emph{multiplicative functions} on $G$,
\item $\CmbA$: the complex of \emph{infinitesimal multiplicative functions} on $A$.
\end{itemize}
We will see several module structures involving these objects, as well as a pair of Van-Est type maps relating them. It will also be explained the relation between these objects and others appearing elsewhere in the literature.

\subsection{The complex $\CmbG$ of multiplicative functions}
\label{sec: CmbG}
Following \cite[\S 9.8]{MackenzieBook05} or \cite[\S 3]{MackenzieX98}, a \emph{multiplicative function} on $G$ is a smooth function $F \in C^{\infty}(G)$
for which $F(gh) = F(g) + F(h)$ for all pairs $(g,h) \in G_2$.
Equivalently, $F : G \to \bR$ is a morphism of Lie groupoids, where the Lie group $\bR$ is considered as a Lie groupoid over a point. Alternatively, $F$ is a 1-cocycle in the complex $\CdG$ computing the differentiable groupoid cohomology $\HdG$ of $G$ \cite[\S 1.2]{Crainic03}. 
Multiplicative functions constitute a subspace $\CmG \subseteq \CG$.
If $f \in \CM$ then $t^*f - s^*f \in \CmG$.

\begin{defi}
\label{def: CmbG}
The 2-term cochain complex $\CmbG$ is defined as follows: $C_m^{0}(G) = C^{\infty}(M)$, $C_m^{1} = C^{\infty}_{m}(G)$, $C_m^{i} = 0$ for $i \ne 0,1$, and the differential $\delta : C^{\infty}(M) \to C_m^{\infty}(G)$ is $f \mapsto t^*f - s^*f$.
The cohomology of $\CmbG$ carries a graded commutative product given by the restriction to $C^{\infty}(M)^G:=\mathrm{Ker} \, \delta $ of the commutative product in $C^{\infty}(M)$ and by $f \cdot [F] := [(t^*f)F]$ for $f \in \CM^G$ and $F \in \CmG$.
\end{defi}

\begin{remi}
The cochain complex $\CmbG$ is equal to the 2-term truncation $C_d^0(G) \to Z_d^1(G)$ of $\CdG$.
In particular, 
\[
H^{0}(\CmbG) \oplus H^{1}(\CmbG) = H^{0}(\CdG) \oplus H^{1}(\CdG)
\] 
as commutative graded algebras.
Note that although $\CdG$ is a dga \cite[\S 1.3]{Crainic03}, the truncated complex $\CmbG$ only carries a graded algebra structure at the level of cohomology. But if $f \in \CM^G$ and $F \in \CmG$, then we have $(t^*f)F \in \CmG$.
\end{remi}

If $H \to G$ is a Morita morphism of Lie groupoids then the induced cochain map $\CdG \to C_d^{\bullet}(H)$ is a quasi-isomorphism \cite[Thm.\ 1]{Crainic03}.
If $G$ is proper then by \cite[Prop.\ 1]{Crainic03}, $H_d^{i}(G) = 0$ for $i \ge 1$, and in particular $H^1(\CmbG) = H_d^1(G) = 0$.

\subsection{The dgla $\XmbG$ of multiplicative vector fields}
\label{sec: XmbG}
Following \cite[\S 9.8]{MackenzieBook05} or \cite[\S 3]{MackenzieX98}, a \emph{multiplicative vector field} is a pair $(X,X_M)$ where $X \in \XG$ and $X_M \in \XM$, which defines a morphism of Lie groupoids $G \to TG$. As a morphism of Lie groupoids $(X,X_M)$ is determined by $X$, and so we often drop $X_M$ from the notation.
This is equivalent to the condition that $X$ is both $s$ and $t$-related to $X_M$, and that the vector field $(X,X)$ on $G_2$ is $\mathrm{pr}_1$, $\mathrm{pr}_2$, and $m$-related to $X$.
The vector space $\XmG$ of multiplicative vector fields is a Lie subalgebra of $\XG$.
If $X \in \XmG$ and $\alpha \in \GA$ then $\alpha^r - \alpha^l \in \XmG$, and $[X,\alpha^r]$ is right invariant so that $[X,\alpha^r] |_M \in \GA$.

\begin{defi}
\label{def: XmbG}
The 2-term dgla $\XmbG$ is defined as follows: $\mathfrak{X}_m^{-1}(G) = \GA$, $\mathfrak{X}_m^{0}(G) = \XmG$, $\mathfrak{X}_m^{i}(G) = 0$ for $i \ne -1,0$, the differential $\p : \GA \to \XmG$ is $\alpha \mapsto \alpha^r - \alpha^l$, and the graded Lie bracket $\LL \cdot, \cdot\RR$ is given by the Lie bracket of multiplicative vector fields and $\LL X,\alpha \RR = [X,\alpha^r] |_{M}$ for $X \in \XmG$ and $\alpha \in \GA$.
\end{defi}

As for any dgla, the cohomology $\HXmbG$ inherits a gla structure from the dgla structure on $\XmbG$ and $H^0(\mathfrak{X}_m(G))$ is a Lie algebra.

\begin{thm}\emph{(\cite{OrtizW19,BerwickEvansL20}.)}
\label{thm: XmbG}
If $H$ is a Lie groupoid Morita equivalent to $G$ then $\XmbG$ and $\XmbH$ are quasi-isomorphic dgla's.
The category 
\[
\GA \ltimes \XmG \rightrightarrows \XmG
\] 
associated to $\XmbG$ is isomorphic to the category $\Gamma(TG)$ of multiplicative sections of the tangent groupoid $TG$ so that $\Gamma(TG)$ carries a natural Lie 2-algebra structure. Also, $\Gamma(TG)$ is equivalent to the category $\Gamma(T\mathcal{X})$ of sections of the tangent stack $T \mathcal{X}$ where $\mathcal{X} = BG$ is the classifying stack of $G$.
\end{thm}

\begin{remi}
The definitions of the categories $\Gamma(TG)$ and $\Gamma(T\mathcal{X})$ can be found in \cite[Def.\ 4.2, 4.3, 4.14]{Hepworth09} and their equivalence as set-theoretic categories is in \cite[Thm.\ 4.15]{Hepworth09}.
Note that in \cite{OrtizW19} and \cite{BerwickEvansL20} the language of Lie algebra crossed modules and of Lie 2-algebras is used instead of dgla's, see \S \ref{sec: 2vec} \& \ref{sec: 2lie} for a dictionary.
A version of Theorem \ref{thm: XmbG} has been generalised to multi-vector fields in \cite[Thm.\ 2.8]{BonechiCLGX20}.
\end{remi}

\begin{remi}
The underlying cochain complex of $\XmbG$ is equal to the 2-term truncation of the \emph{deformation complex} $\CdefG$ of $G$ shifted by 1, defined in \cite[Def.\ 2.1]{CrainicMS20}. The deformation cohomology of $G$ is denoted $\HdefG$.
In particular, 
\[
H^{-1}(\XmbG) \oplus H^{0}(\XmbG) = H^{0}_{\mathrm{def}}(G) \oplus H^{1}_{\mathrm{def}}(G)
\] 
as graded vector spaces.
It is shown in \cite[\S 9]{CrainicMS20} that $\HdefG$ is isomorphic to the cohomology of $G$ with coefficients in the adjoint representation up to homotopy of $G$ defined in \cite{AbadC13}. Additionally, if $G$ is proper then $H^i_{\mathrm{def}}(G)=0$ for every $i\geq 2$, \cite[Thm.\ 6.1]{CrainicMS20}.
\end{remi}

\subsubsection{The left $\CM^G$-module structure on $\CmbG$ and $\HCmbX$}
\label{sec: CMG_mod}

We will show now that the complex $\XmbG$ of vector spaces is actually a complex of  $\CMG$-modules.

\begin{prop}
\label{prop: CMG_chcomplex}
Suppose that $f \in \CMG$, $\alpha \in \GA$ and $X \in \XmG$. Then the following hold:
\begin{enumerate}
\item $(t^*f)X$ is multiplicative with $((t^*f)X)_M = fX_M$.
\item $\p(f\alpha) = (t^*f) \p(\alpha)$.
\end{enumerate}
\end{prop}

\begin{proof}
We start by showing \emph{(a)}. It is clear that $(t^*f)X = (s^*f)X$ is both $s$ and $t$-related to $fX_M$.
If $(g,h) \in G_2$, with $s(g) = t(h) = x \in M$, then
\begin{align*}
dm ( ((s^*f)X)_g , ((t^*f)X)_h ) & = f(x) dm (X_g,X_h) \\
& = f(x) X_{gh} \\
& = ((t^*f)X)_{gh} .
\end{align*}

In order to prove \emph{(b)}, we calculate
\begin{align*}
\p (f\alpha) & = (f\alpha)^r - (f\alpha)^l \\
& = (t^*f) \alpha^r - (s^*f) \alpha^l \\
& = (t^*f) \left( \alpha^r - \alpha^l \right) \\
& = (t^*f) \p \alpha,
\end{align*}
finishing the proof.
\end{proof}

As a consequence, one has the following module structures.

\begin{defi}
\label{def: CMG_mod}
The chain complex $\XmbG$ of vector spaces is a chain complex of left $\CMG$-modules with action defined by
\begin{align*}
\CM^G \otimes \XmbG & \to \XmbG \\
f \otimes \left( \alpha,X \right) & \mapsto \left( f\alpha , (t^*f)X \right) .
\end{align*}
This action descends to cohomology to define a graded $\CMG$-module structure on $\HXmbG$ given by
\begin{align*}
\CM^G \otimes \HXmbG & \to \HXmbG \\
f \otimes \left( [\alpha],[X] \right) & \mapsto \left( [f\alpha] , [(t^*f)X] \right) .
\end{align*}
\end{defi}

\subsubsection{The graded module structure of $\HXmbG$ over $\HCmbG$}
\label{sec: HCmbG_mod}
The $\CM^G$-module structure on $\HXmbG$ defined in Definition \ref{def: CMG_mod} actually extends to an $\HCmbG$-module structure.
In \cite[Lem.\ 2.5]{CrainicMS20} it is shown that the deformation complex $\CdefG$ carries a natural right dg-module structure over the dga $\CdG$.
This descends to cohomology to define a right graded $\HdG$-module structure on $\HdefG$.
Truncating this module structure gives the module structure in the following definition.

\begin{defi}
\label{def: cdot}
The left graded module structure of $\HXmbG$ over $\HCmbG$
\begin{align*}
\HCmbG \otimes \HXmbG & \to \HXmbG \\
\left( [f],[F] \right) \otimes \left( [\alpha],[X] \right) & \mapsto \left( [f] \cdot [\alpha], [f] \cdot [X]+[F]\cdot [\alpha] \right) \\
& = \left( [f \alpha] , [(t^*f) X]+[F \alpha^r] \right)
\end{align*}
is defined as follows:
\begin{enumerate}
\item \label{it: def_cdot_FX}
$[F] \cdot [X] = 0$

\item \label{it: def_cdot_Fa}
$[F] \cdot [\alpha] = [F \alpha^r]$

\item \label{it: def_cdot_fX}
$[f] \cdot [X] = [(t^* f)X]$

\item \label{it: def_cdot_fa}
$[f] \cdot [\alpha] = [f\alpha]$
\end{enumerate}
for $F \in \CmG$, $f \in \CM^G$, $X \in \XmG$ and $\alpha \in \mathrm{Ker} \, \p \subseteq \GA$.
\end{defi}

\begin{remi}
Note that as the graded algebra $\HCmbG$ is strictly commutative (as opposed graded commutative) the \textbf{right} $\HdG$-module structure on $\HdefG$ defines a \textbf{left} $\HCmbG$-structure on $\HXmbG$ as in Definition \ref{def: cdot}.
\end{remi}

\begin{remi}\label{rem: Far}
Note that the vector field $F\alpha^r\in \mathfrak{X}(G)$ in Definition \ref{def: cdot}\eqref{it: def_cdot_Fa} is indeed multiplicative. In fact, $\alpha\in \mathrm{Ker}(\partial)$ is equivalent to saying that $\alpha^r=\alpha^l$ and hence $\alpha^r$ is both $s$ and $t$-projectable to zero. It remains to check that $(F\alpha^r,F\alpha^r)$ is $m$-related to $F\alpha^r$. For that we use the explicit formula of the multiplication in $TG$ as in \cite[Thm 1.4.14]{MackenzieBook05}. Given composable arrows $g,h$ and local bisections $\sigma:U\to G$ and $\tau:V\to G$ with $s(g)\in U$ and $s(h)\in V$, then

$$(F\alpha^r)_g\bullet (F\alpha^r)_h=dL_{\sigma}((F\alpha^r)_h)+dR_{\tau}((F\alpha^r)_g)-dL_{\sigma}dR_{\tau}d1ds(F\alpha^r)_g.$$ 

\noindent The last term of the right hand side vanishes since $\alpha^r$ is $s$-projectable to zero. Using that $\alpha^r=\alpha^l$ one easily check that 

$$dL_{\sigma}((F\alpha^r)_h)+dR_{\tau}((F\alpha^r)_g)=F(h)\alpha^l_{gh}+F(g)\alpha^r_{gh}=F(gh)\alpha^r_{gh},$$
where in the last identity we have used that $F$ is multiplicative.

\end{remi}

\subsection{The complex $\CmbA$ of infinitesimal multiplicative functions}
\label{sec: CmbA}
We define the vector space $C_m^{1}(A) \subseteq \Gamma(A^*)$ to be the subspace with elements $\omega$ satisfying 
\begin{equation}
\label{eqn: w_cocycle}
\mL_{a(\alpha)} \omega(\beta) - \mL_{a(\beta)} \omega(\alpha) - \omega([\alpha,\beta]) = 0
\end{equation}
for all $\alpha,\beta \in \GA$. In other words $\omega : A \to \bR$ is a morphism of Lie algebroids where $\bR$ is considered as an abelian Lie algebra.
We think of elements of $C_m^1(A)$ as ``infinitesimal multiplicative functions'' on $A$.
If $f \in \CM$ then $d_Af\in C_m^{1}(A)$ is defined as $d_Af(\alpha) :=  \mL_{a(\alpha)}f$.

\begin{defi}
\label{def: CmbA}
The 2-term complex of infinitesimal multiplicative functions $\CmbA$ is defined as follows: $C_m^{0}(A) = C^{\infty}(M)$, $C_m^{1}(A)$ is as defined above, $C_m^{i}(A) = 0$ for $i \ne 0,1$, and the differential is $d_A : C^{\infty}(M) \to C_m^1(A)$. 
The cohomology of $\CmbA$ carries a graded product given by $[f][f'] = [ff']$ and $[f] [\omega] := [f\omega]$ for $f,f' \in \CM$ and $\omega \in C_m^1(A)$.
\end{defi}

\begin{remi}
The complex $\CmbA$ is equal to the 2-term truncation $C^0(A) \to Z^1(A)$ of the 
\emph{Chevalley-Eilenberg complex} $\CA$ of $A$, the cohomology of which is denoted $\HA$, see \cite[\S 7]{MackenzieBook05} or \cite[\S 1.4]{Crainic03}.
In particular,
\[
H^{0}(\CmbA) \oplus H^{1}(\CmbA) = H^{0}(A) \oplus H^{1}(A)
\]
as graded vector spaces.
Note that although $\CA$ carries a natural dga structure, $\CmbA$ only carries a graded algebra structure at the level of cohomology.
\end{remi}

\subsection{The dgla $\XmbA$ of derivations of $A$}
\label{sec: XmbA}
A \emph{Lie algebroid derivation of} $A$ is a pair $(D,\sigma(D))$ where $D \in \mathrm{Der}(\GA)$ is a derivation of the Lie algebra $\GA$ and $\sigma(D) \in \XM$ with $D(f\alpha) = (\mL_{\sigma(D)}f) \alpha + fD(\alpha)$ and $\rho(D\alpha)=[\sigma(D),\rho(\alpha)]$ for all $f \in \CM$ and $\alpha \in \GA$ \cite[Def.\ 4.1]{MoerdijkM02}. The vector field $\sigma(D)$, called the \emph{symbol} of $D$, is uniquely determined by $D$ and so we sometimes drop it from the notation.
The space of Lie algebroid derivations is a Lie subalgebra $\mathrm{Der}(A) \subseteq \mathrm{Der}(\GA)$ and makes part of a 2-term dgla as explained below.

\begin{defi}
\label{def: XmbA}
The 2-term dgla $\XmbA$ is defined as follows: $\mathfrak{X}_m^{-1}(A) = \GA$, $\mathfrak{X}_m^{0} = \mathrm{Der}(A)$, $\mathfrak{X}_m^{i} = 0$ for $i \ne -1,0$, the differential $\mathrm{ad} : \GA \to \mathrm{Der}(A)$ is $\alpha \mapsto ([\alpha,-] , a(\alpha))$, and the graded Lie bracket is given the commutator of Lie algebroid derivations and $[(D,\sigma(D)) , \alpha] = D(\alpha)$
for $(D,\sigma(D)) \in \mathrm{Der}(A)$ and $\alpha \in \GA$.
\end{defi}

\begin{remi}
The dgla $\XmbA$ is the 2-term truncation of the \emph{deformation complex} $\CdefA$ of $A$ shifted by 1, defined in \cite[\S 2]{CrainicM08}. The deformation cohomology of $A$ is denoted $\HdefA$. 
The dgla structure on $\XmbA$ is a restriction of the dgla structure defined in loc.\ cit.\ on the shift by 1 of $\CdefA$.
In particular, 
\[
H^{-1}(\XmbA) \oplus H^{0}(\XmbA) = H^{0}_{\mathrm{def}}(A) \oplus H^{1}_{\mathrm{def}}(A)
\]
as graded vector spaces.
\end{remi}

\subsubsection{The dgla $\XmbA$-module structure on $\CmbA$.}
The deformation complex $\CdefA$ is isomorphic to the dgla of graded derivations of the dga $\CA$ and so $\CA$ carries a natural $\CdefA$-module structure \cite[\S 2.5-4.8]{CrainicM08}.
This module structure is compatible with the inclusions $\XmbA \subseteq \CdefA$ and $\CmbA \subseteq \CA$, which gives rise to the following.

\begin{defi}
\label{def: obull}
The dgla module structure of $C_{m}^{\bullet}(A)$ over $\mathfrak{X}_{m}^{\bullet}(A)$ is defined as:
\begin{align*}
\overline{\mu} : \XmbA \otimes \CmbA & \to \CmbA \\
(\alpha , D) \otimes (f,\omega) & \mapsto (\alpha \obullet \omega + D \obullet f , D \obullet \omega)
\end{align*}
where
\begin{enumerate}
\item \label{it: def_obull_Dw} 
$D \obullet \omega  = \mL_{\sigma(D)} \circ \omega - \omega \circ D$
\item \label{it: def_obull_Df} 
$D \obullet f  = \mL_{\sigma(D)} f$
\item \label{it: def_obull_aw} 
$\alpha \obullet \omega  = i_{\alpha} \omega := \omega(\alpha)$
\item \label{it: def_obull_af} 
$\alpha \obullet f  = 0$
\end{enumerate}
for $D \in \mathrm{Der}(A)$, $\alpha \in \GA$, $\omega \in C_m^{1}(A)$ and $f \in \CM$.
\end{defi}

In Theorem \ref{thm: mu} we prove that there is a related dgla module structure of $C_{m}^{\bullet}(G)$ over $\XmbG$.

\subsubsection{The graded module structure of $\HXmbA$ over $H^{\bullet}(C_{m}^{\bullet}(A))$.}
The deformation complex $\CdefA$ of $A$ is isomorphic to the complex $\Omega^{\bullet}(A,\mathrm{ad})$ computing the cohomology of the adjoint representation up to homotopy of $A$ \cite[Thm.\ 3.11]{AbadC12}.
The latter complex is, by definition \cite[Def.\ 3.1]{AbadC12}, a module over the differential graded algebra $\CA$.
Via the isomorphism $\CdefA \cong \Omega^{\bullet}(A,\mathrm{ad})$ this module structure descends to cohomology to define a right graded $\HA$-module structure on $\HdefA$, and then applying appropriate truncations gives the module structure in the following definition.

\begin{defi}
\label{def: ocdot}
The left graded $H^{\bullet}(C_{m}^{\bullet}(A))$-module structure on $H^{\bullet}(\mathfrak{X}_{m}^{\bullet}(A))$
\begin{align*}
\HCmbA \otimes \HXmbA & \to \HXmbA \\
\left( [f],[\omega] \right) \otimes \left( [\alpha],[D] \right) & \mapsto \left( [f] \ocdot [\alpha], [\omega] \ocdot [\alpha] + [f] \ocdot [D] \right)
\end{align*}
is defined as follows:
\begin{enumerate}
\item \label{it: def_ocdot_FD}
$[\omega] \ocdot [D] = 0$

\item \label{it: def_ocdot_wa}
$[\omega] \ocdot [\alpha] = [\omega \alpha]$, where $\omega \alpha := (\beta \mapsto \omega(\beta) \, \alpha)$

\item \label{it: def_ocdot_fD}
$[f] \ocdot [D] = [fD]$

\item \label{it: def_ocdot_fa}
$[f] \ocdot [\alpha] = [f\alpha]$
\end{enumerate}
for $\omega \in C_m^1(A)$, $f \in \mathrm{Ker}\,d_A = H^{0}(\CmbA)$, $D \in \mathrm{Der}(A)$ and $\alpha \in \mathrm{Ker} \, \mathrm{ad} = H^{-1}(\XmbA)$.
\end{defi}

\begin{remi}
\label{rem: wa}
Just as in Remark \ref{rem: Far},
one observes that $\omega \alpha$ in Definition \ref{def: ocdot}\eqref{it: def_ocdot_wa} is indeed a Lie algebroid derivation of $A$.
If $\beta,\beta' \in \GA$ then using the fact that $\alpha \in \mathrm{Ker} \, \p$ implies that $\alpha$ is central in $\GA$, and \eqref{eqn: w_cocycle}:
\begin{align*}
[ (\omega \alpha) (\beta), \beta'] + [\beta, (\omega \alpha) (\beta')] 
& = [ \omega(\beta) \alpha , \beta'] + [\beta, \omega(\beta') \alpha] \\
& = \omega(\beta)[\alpha,\beta'] + \omega(\beta')[\beta,\alpha] + (\mL_{a(\beta)} \left(\omega(\beta')) - \mL_{a(\beta')} (\omega(\beta)) \right) \alpha \\
& = \omega([\beta,\beta']) \alpha \\
& = (\omega \alpha) ([\beta,\beta']),
\end{align*}
which shows that $\omega \alpha$ is a derivation of the Lie algebra $\GA$. By Definition, $\omega \alpha$ is $\CM$-linear, and so is a Lie algebroid derivation with zero symbol.
\end{remi}

\subsection{Van-Est maps}
\label{sec: van_est}
If $F \in \CmG$ then the function $\omega_F : \alpha \mapsto u^* \mL_{\alpha^r}(F)$ is an element of $C_m^1(A)$.
This is just the Lie functor mapping a Lie groupoid morphism $F : G \to \bR$ to the corresponding Lie algebroid morphism $A \to \bR$, where in the first case $\bR$ is considered as a Lie group and in the second as its Lie algebra.

\begin{defi}
\label{def: oVE}
The \emph{Van-Est map} $\oVE : \CmbG \to \CmbA$ is the morphism of cochain complexes that is the identity on $\CM$ and is $C_m^1(G) \to C_m^1(A)$, $F \mapsto \omega_F$ in degree 1.
\end{defi}

By Lie's theorems for Lie groupoids, $\oVE$ is injective (resp.\ an isomorphism) if $G$ is source connected (resp.\ source simply connected).
The map $\oVE$ is the truncation of the Van-Est map $C_d^{\bullet}(G) \to C^{\bullet}(A)$ defined in \cite{WeinsteinX91}, whose associated map on cohomology $H^i(G) \to H^i(A)$ is an isomorphism in degrees $i \le n$ and injective in degree $i=n+1$ whenever $G$ is source $n$-connected \cite[Thm.\ 4]{Crainic03}.\\

If $X \in \XmG$ then the map $D_X : \alpha \mapsto [X,\alpha^r]|_M$ is a derivation of $\GA$ and the pair $(D_X,X_M)$ is an element of $\mathrm{Der}(A)$.
This defines a Lie algebra homomorphism from $\XmG$ to $\mathrm{Der}(A)$ \cite[\S 3]{MackenzieX98}, \cite[\S 4]{MoerdijkM02}.
\begin{defi}
\label{def: VE}
The \emph{Van-Est map} $\VE : \XmbG \to \XmbA$ is the morphism of cochain complexes that is the identity map on $\GA$ and is $\XmG \to \mathrm{Der}(A)$, $X \mapsto D_X$ in degree 0. 
\end{defi}

The Lie algebra morphism $X \mapsto D_X$ and therefore the morphism $\VE : \XmbG \to \XmbA$ is injective (resp.\ an isomorphism) whenever $G$ is source-connected (resp.\ source simply connected) \cite[Thm.\ 4.5]{MoerdijkM02}. Also, up to a shift in degree, the map $\VE$ is the truncation of the Van-Est map $\CdefG \to \CdefA$ defined in \cite[\S 10]{CrainicMS20} whose associated map on cohomology $H_{\mathrm{def}}^i(G) \to H_{\mathrm{def}}^i(A)$ is an isomorphism in degrees $i \le n-1$ whenever $G$ is source $n$-connected \cite[Thm.\ 10.1]{CrainicMS20}, \cite[Thm.\ 4.7]{AbadS11}.

One easily observes that $\VE$ is a morphism of dgla's. This follows from a standard construction for general dgla's as explained in Remark \ref{rem: adrep} below.

\begin{remi}
\label{rem: adrep}
Suppose that $L^{\bullet}=L^{-1} \oplus L^{0}$ is a 2-term dgla with differential $d_L$.
Then $L^{-1}$ carries a Lie bracket $[-,-]_{L^{-1}}$ defined by $[x,x']_{L^{-1}} := [d_{L} x,x']$, and $L^{0}$ acts on this Lie algebra by derivations via $y \mapsto [y,-] \in \mathrm{Der}(L^{-1})$.
These operations define a second dgla $L^{-1} \oplus \mathrm{Der}(L^{-1})$ with differential $x \mapsto \mathrm{ad}_x$, and a morphism of dgla's $L^{-1} \oplus L^{0} \to L^{-1} \oplus \mathrm{Der}(L^{-1})$.

In the case of the dgla $\XmbG$, the map $\XmG \to \mathrm{Der}(\GA)$, $X \mapsto [X,-]$ lands in the Lie subalgebra $\mathrm{Der}(A) \subseteq \mathrm{Der}(\GA)$ and so one gets a morphism of dgla's $\XmbA \to \XmbA$, which is exactly the map $\VE$ of Definition \ref{def: VE}.
\end{remi}

\section{Functions on differentiable stacks}

There are several equivalent descriptions of the algebra $\CX$ of smooth functions on a differentiable stack $\mathcal{X}$. These are described in the following proposition.

\begin{prop}
\label{prop: fun}
Let $\mathcal{X}$ be a differentiable stack, $M \to \mathcal{X}$ an atlas, and $G$ the associated Lie groupoid over $M$.
The following commutative algebras are canonically isomorphic:
\begin{enumerate}
\item morphisms of stacks $\mathrm{Hom}_{\St} (\mathcal{X} , \bR)$ from $\mathcal{X}$ to the manifold $\bR$.
\item global sections $H^{0}(C^{\infty}_{\mathcal{X}})$ of the sheaf of smooth functions on $\mathcal{X}$.
\item invariant functions $\CMG$ on $M$.
\item the degree zero differentiable cohomology $H_d^{0}(G)$  of $G$.
\end{enumerate}
\end{prop}

\begin{proof}
See \cite[\S 3]{BehrendX03}.
\end{proof}

The categories defined below are directly related to the previous algebras.

\begin{thm}
\label{thm: fun}
Let $G$ be a Lie groupoid over $M$.
The following two categories are isomorphic:
\begin{enumerate}
\item \label{it: thm_fun_CmbG}
The category $\CM \ltimes \CmG$ associated to the 2-term complex $\CM \xrightarrow{\delta} \CmG$.

\item \label{it: thm_fun_HomGpd} 
The category $\mathrm{Hom}_{\Gpd} (G,\bR)$ of Lie groupoid homomorphisms from $G$ to $\bR$ considered as a Lie groupoid over a point.
\end{enumerate}
In addition, if $M$ is Hausdorff and paracompact then \eqref{it: thm_fun_CmbG} and \eqref{it: thm_fun_HomGpd} are equivalent to:
\begin{enumerate}
\item[(c)] \label{it: thm_fun_HomSt}
The category $\mathrm{Hom}_{\St} (BG,B\bR)$ of morphisms of stacks from $BG$ to $B\bR$.
\end{enumerate}
More generally, if $\mathcal{X}$ is a differentible stack, $M \to \mathcal{X}$ is an atlas with $M$ Hausdorff and paracompact, and $G$ is the associated Lie groupoid over $M$, then the categories \eqref{it: thm_fun_CmbG} and \eqref{it: thm_fun_HomGpd} are equivalent to $\mathrm{Hom}_{\St} (\mathcal{X},B\bR)$. 
\end{thm}

\begin{proof}
We first show that the categories \eqref{it: thm_fun_CmbG} and \eqref{it: thm_fun_HomGpd} are isomorphic.
This can be deduced from \cite[Thm.\ 3.1]{OrtizW19}, see Remark \ref{rem: VBgpds} below, but we give a self-contained proof.

The set of objects of the category $\mathrm{Hom}_{\Gpd} (G,\bR)$ is equal to the set $\CmG$ of multiplicative functions on $G$ (see \S \ref{sec: CmbG}).
A morphism $f : F \Rightarrow F'$ in the category $\mathrm{Hom}_{\Gpd} (G,\bR)$ is smooth natural transformation; that is, a smooth function $f : M \to \bR$ such that for each $g \in G$, with $s(g)=x$ and $t(g)=y$, the following square commutes:
\[
\xymatrix{
F(x) \ar[r]^{F(g)} \ar[d]_{f(x)} & F(y) \ar[d]^{f(y)} \\
F'(x) \ar[r]_{F'(g)} & F'(y)
}
\]
which is equivalent $F'(g) = F(g) + f(y) - f(x)$.
It follows that a morphism $f : F \Rightarrow F'$ is exactly a smooth function $f \in \CM$ such that $F' = F + t^*f - s^*f$, or $F' = F + \delta f$.
This shows that the map $(F',f,F) \mapsto f : F \Rightarrow F'$ is a bijection from the set of morphisms in $\mathrm{Hom}_{\Gpd} (G,\bR)$ to the set of morphisms in $\CM \ltimes \CmG$.

The composition of smooth natural transformations $f: F \Rightarrow F'$ and $f' : F' \Rightarrow F''$ is given by $(f \circ f') (x) = f(x) + f'(x)$, which corresponds exactly to the composition $(F'',f',F') \circ (F',f,F) = (F'',f'+f,F)$ in the category $\CM \ltimes \CmG$.\\
\\
It remains to show that the categories \eqref{it: thm_fun_CmbG} and \eqref{it: thm_fun_HomGpd} are equivalent to $\mathrm{Hom}_{\St} (BG,B\bR)$.
For any Lie groupoid $H$ there is a natural equivalence of categories 
\begin{equation}
\label{eqn: GPD_St}
\mathrm{Hom}_{\GPD} (G,H) \xrightarrow{\simeq} \mathrm{Hom}_{\St} (BG,BH) .
\end{equation}
where $\mathrm{Hom}_{\GPD} (G,H)$ is the category of $G$-$H$ bibundles or Hilsum-Skandalis morphisms, and a natural functor
\begin{equation}
\label{eqn: Gpd_GPD}
\mathrm{Hom}_{\Gpd} (G,H) \to \mathrm{Hom}_{\GPD} (G,H) 
\end{equation}
which maps a morphism of Lie groupoids to the corresponding $G$-$H$ bibundle; see \cite[\S 3.2-3.3]{Lerman10} for further details.
The functor \eqref{eqn: Gpd_GPD} is fully-faithful because natural transformations between morphisms of Lie groupoids correspond exactly to isomorphisms between the corresponding bibundles, and by \cite[Lem.\ 3.36]{Lerman10} the essential image of \eqref{eqn: Gpd_GPD} is the full subcategory of $\mathrm{Hom}_{\GPD} (G,H)$ consisting of bibundles that admit a global section.

Combining these facts, it is sufficient to show that if $P$ is a $G$-$\bR$ bibundle then the principal $\bR$-bundle $P \to M$ admits a global section.
This follows from the assumptions on the topology of the manifold $M$: being Hausdorff and paracompact implies that the sheaf $C_M^{\infty}$ of smooth $\bR$-valued functions on $M$ is fine and therefore the \v{C}ech cohomololgy group $\check{H}^1(M,C_M^{\infty})$ classifying principal $\bR$-bundles is zero.
\end{proof}

\begin{remi}
\label{rem: VBgpds}
The isomorphism $\CM \ltimes \CmG \simeq \mathrm{Hom}_{\Gpd} (G,\bR)$ can also be deduced from results in \cite{OrtizW19}.
In the category $\Gpd$, morphisms from $G$ to $\bR$ are equivalent to sections of the projection $G \times \bR \to G$, which is a $\mathcal{VB}$-groupoid over $G$ with core the trivial vector bundle $M \times \bR$, see loc.\ cit.\ for the terminology.
The isomorphism then follows from \cite[Thm.\ 3.1]{OrtizW19} after identifying the right (resp.\ left) invariant function $f^r$ (resp.\ $f^l$) on $G$ associated to $f \in \Gamma(M \times \bR) = \CM$ with $t^* f$ (resp.\ $s^*f$).
\end{remi}

\section{Some identities}
\label{sec: ident}
Throughout \S\ref{sec: ident} we use the notation established in \S \ref{sec:bg_notation}, and we use $f,f'$ to denote elements of $\CM$; $F,F'$ to denote elements of $\CmG$; $\alpha,\beta$ to denote elements of $\GA$; and $X,Y$ to denote elements of $\XmG$.

\subsection{The operation `$\bullet$'}
In order to prove the main results of \S\ref{sec: mod} and \S\ref{sec: co} we first establish a number of relevant identities satisfied by the following operations.

\begin{defi}
\label{def: bull}
We define a canonical degree zero map between dgla's  \[\mathfrak{X}_m^{\bullet}(G)\to \mathrm{End}^{\bullet}(C^{\bullet}_m(G)),\] given by:
\begin{enumerate}
\item \label{it: def_bull_Xf} 
$X \bullet f  = \mL_{X_M}f$;
\item \label{it: def_bull_XF} 
$X \bullet F  = \mL_{X}F$;
\item \label{it: def_bull_af} 
$\alpha \bullet f  = 0$;
\item \label{it: def_bull_aF} 
$\alpha \bullet F  = u^*\mL_{\alpha^r}F$;
\end{enumerate}
for every $X \in \mathfrak{X}_m(G)$, $\alpha \in \Gamma(A)$, $f \in C^{\infty}(M)$ and $F \in C^{\infty}_m(G)$.
\end{defi}

\begin{remi}
Note that $\mL_{X}F$ is a multiplicative function whenever $X$ and $F$ are multiplicative because \[m^*(\mathcal{L}_XF)=\mathcal{L}_{(X,X)}(m^*F)=\mathcal{L}_{(X,X)}(\mathrm{pr}_1^*F+\mathrm{pr}_2^*F)=\mathrm{pr}_1^*(\mathcal{L}_{X}F)+\mathrm{pr}_2^*(\mathcal{L}_{X}F).\]
\end{remi}

In other words, associated to each multiplicative vector field (resp.\ section of $A$) there is a degree zero (resp.\ minus one) endomorphism of the graded vector space $C^{\infty}(M) \oplus C^{\infty}_m(G)$ underlying the complex $\CmbG$ defined as above.

\subsection{Derivatives of multiplicative functions}
We need the following identities about derivatives of multiplicative functions to prove several results about the operation `$\bullet$' of Definition \ref{def: bull}.

\begin{lem}
\label{lem: aF}
The following identities hold:
\begin{enumerate}
\item \label{it: lem_aF_rl}
$\mL_{\alpha^l} F = i^*(\mL_{\alpha^r} F)$.

\item \label{it: lem_aF_rinv}
$(\mL_{\alpha^r} F)(h) = (\mL_{\alpha^r} F)(hg)$ whenever $g,h \in G$ with $s(h)=t(g)$.

\item \label{it: lem_aF_linv}
$(\mL_{\alpha^l} F)(h) = (\mL_{\alpha^l} F)(gh)$ whenever $g,h \in G$ with $t(h)=s(g)$.

\item \label{it: lem_aF_tu}
$\mL_{\alpha^r} F = t^*u^* \mL_{\alpha^r} F = t^*u^*\mL_{\alpha^l}F$.

\item \label{it: lem_aF_su}
$\mL_{\alpha^l} F = s^*u^* \mL_{\alpha^l} F = s^*u^* \mL_{\alpha^l} F$.
\end{enumerate}
\end{lem}

\begin{proof}
\eqref{it: lem_aF_rl}.
Using the fact that $i^*F = -F$ and that $i_* \alpha^r = -\alpha^l$, where $i_*\alpha^r$ is the pushforward of $\alpha^r$ by the diffeomorphism $i$, we have that
\begin{align*}
\mL_{\alpha^r} F & = - \mL_{\alpha^r} i^* F \\
& = - i^* (\mL_{i_* \alpha^r} F) \\
& = i^*(\mL_{\alpha^l} F).
\end{align*}

\eqref{it: lem_aF_rinv}.
Suppose that $g \in G$ with $x=s(g)$ and $y=t(g)$.
First note that as $\alpha^r$ is tangent to the fibers of the submersion $s: G \to M$ it restricts to a vector field $\alpha^r|_{s^{-1}(z)}$ on the submanifold $s^{-1}(z)$ for each $z \in M$, and the right invariance of $\alpha^r$ is equivalent to the condition that $(R_g)_* (\alpha^r|_{s^{-1}(y)}) = \alpha^r|_{s^{-1}(x)}$.
Moreover, the value of $\mL_{\alpha^r} F$ at $h \in G$ with $s(h)=z$ only depends on the restriction of $F$ to $s^{-1}(z)$, that is 
\[
(\mL_{\alpha^r}F)(h) = \left( \mL_{\alpha^r|_{s^{-1}(z)}} \left(F|_{s^{-1}(z)}\right) \right) (h) .
\]

Next, if follows from the multiplicativity of $F$ that $F|_{s^{-1}(y)}$ is equal to $R_g^* (F|_{s^{-1}(x)})$ up to the addition of a constant function:
\begin{align*}
R_g^* (F|_{s^{-1}(x)}) (h) & = F|_{s^{-1}(x)} (hg) \\
& = F(h) + F(g) \\
& = F|_{s^{-1}(y)}(h) + F(g).
\end{align*}
This implies that $\mL_{Z} R_g^* (F|_{s^{-1}(y)}) = \mL_{Z} (F|_{s^{-1}(y)})$ for any vector field $Z$ on $s^{-1}(y)$.

Combining these two observations and using the fact that $\alpha^r$ is right invariant we have that
\begin{align*}
(\mL_{\alpha^r} F) |_{s^{-1}(y)} & = \mL_{\alpha^r|_{s^{-1}(y)}} (F|_{s^{-1}(y)}) \\
& = \mL_{\alpha^r|_{s^{-1}(y)}} (R_g^* (F|_{s^{-1}(x)}) ) \\
& = R_g^* \left( \mL_{(R_g)_* (\alpha^r|_{s^{-1}(y)}) } F|_{s^{-1}(x)} \right) \\
& = R_g^* \left( \mL_{\alpha^r|_{s^{-1}(x)}} F|_{s^{-1}(x)} \right)
\end{align*}
and therefore $(\mL_{\alpha^r}F) (h) = (\mL_{\alpha^r}F) (hg)$ whenever $h \in G$ with $s(h)=t(g)$.

\eqref{it: lem_aF_linv}.
This follows from the same argument as in the proof of \eqref{it: lem_aF_rinv}, but with $\alpha^r$ replaced by $\alpha^l$ and $R_g$ replaced by $L_g$.

\eqref{it: lem_aF_tu}.
Using part \eqref{it: lem_aF_rinv} we have that
\begin{align*}
(\mL_{\alpha^r}F)(g) & = (\mL_{\alpha^r}F)(1_{t(g)}g) \\
& = (\mL_{\alpha^r}F)(1_{t(g)}) \\
& = (t^*u^*\mL_{\alpha^r}F)(g)
\end{align*}
which proves the first equality.
To prove the second equality we use part \eqref{it: lem_aF_rl} and the fact that $iu=u$:
\begin{align*}
t^*u^*\mL_{\alpha^r}F & = t^*u^*i^*\mL_{\alpha^l}F \\
& = t^*(iu)^*\mL_{\alpha^l}F \\
& = t^*u^*\mL_{\alpha^l}F.
\end{align*}

\eqref{it: lem_aF_su}.
This follows from part \eqref{it: lem_aF_linv} and the same argument as in the proof of \eqref{it: lem_aF_tu}.
\end{proof}

Lemma \ref{lem: aF} \eqref{it: lem_aF_tu}\&\eqref{it: lem_aF_su} can be restated in terms of Definition \ref{def: bull} as follows.

\begin{lem}[Derivatives]
\label{lem: LaF}
The following statements hold:
\begin{enumerate}
\item \label{it: lem_LaF_r} 
$\mL_{\alpha^r} F = t^* (\alpha \bullet F)$

\item \label{it: lem_LaF_s} 
$\mL_{\alpha^l} F = s^* (\alpha \bullet F)$
\end{enumerate}
\end{lem}

Similarly, the compatibility between the operation $\bullet$ in Definition \ref{def: bull} and the involved differential can be expressed as follows.

\begin{lem}[Differentials]
\label{lem: diff}
The following identities hold:
\begin{enumerate}
\item \label{it: lem_diff_dXf}
$\delta (X \bullet f) = X \bullet \delta f$

\item \label{it: lem_diff_adf}
$\alpha \bullet \delta f  = \p \alpha \bullet f$

\item \label{it: lem_diff_daF}
$\delta (\alpha \bullet F)  = \p \alpha \bullet F$
\end{enumerate}
\end{lem}

\begin{proof}
\eqref{it: lem_diff_dXf}. This follows from the fact that $X$ is $s$ and $t$-related to $X_M$.

\eqref{it: lem_diff_adf}. Using the fact that $\alpha^r$ is $s$-related to zero and $t$-related to $a(\alpha)$, and that $tu = \mathrm{id}_{M}$, we have that
\begin{align*}
\alpha \bullet \delta f & = u^* \mL_{\alpha^r} t^*f - u^* \mL_{\alpha^r} s^*f \\
& = u^* t^* \mL_{a(\alpha)} f \\
& = (tu)^* \mL_{a(\alpha)} f \\
& = \mL_{a(\alpha)} f \\
& = \p \alpha \bullet f.
\end{align*}

\eqref{it: lem_diff_daF}. Using Lemma \ref{lem: aF}\eqref{it: lem_aF_tu} and \eqref{it: lem_aF_su} we have:
\begin{align*}
\delta(\alpha \bullet F) & = t^*u^*\mL_{\alpha^r}F - s^*u^*\mL_{\alpha^r}F \\
& = \mL_{\alpha^r}F - \mL_{\alpha^l}F \\
& = \mL_{\alpha^r - \alpha^l}F \\
& = \mL_{\p(\alpha)}F \\
& = \p \alpha \bullet F.
\qedhere
\end{align*}
\end{proof}

Finally, recall the Lie structure on $\XmbG$ from Definition \ref{def: XmbG}, that $\LL X,Y \RR = [X,Y]$ and $\LL X,\alpha \RR = [X,\alpha^r]|_M$. Hence the operation $\bullet$ in Definition \ref{def: bull} is compatible with Lie brackets in the following sense.

\begin{lem}[Lie brackets]
\label{lem: lie}
The following identities hold:
\begin{enumerate}
\item \label{it: lem_lie_XYf}
$\LL X,Y \RR  \bullet f  = X \bullet (Y \bullet f) - Y \bullet (X \bullet f)$

\item \label{it: lem_lie_XYF}
$\LL X,Y \RR  \bullet F  = X \bullet (Y \bullet F) - Y \bullet (X \bullet F)$

\item \label{it: lem_lie_XaF}
$\LL X,\alpha \RR \bullet F  = X \bullet (\alpha \bullet F) - \alpha \bullet (X \bullet F)$.
\end{enumerate}
\end{lem}

\begin{proof}
Statements \eqref{it: lem_lie_XYf} and \eqref{it: lem_lie_XYF} follow immediately from standard properties of Lie derivatives and the fact that $\LL X,Y \RR = [X,Y]$ and $\LL X,Y \RR _M = [X_M,Y_M]$.

\eqref{it: lem_lie_XaF}. This follows from standard properties of Lie derivatives and the fact that $X$ is $u$-related to $X_M$:
\begin{align*}
\LL X,\alpha \RR \bullet F - X \bullet (\alpha \bullet F) + \alpha \bullet (X \bullet F) 
& = (u^*\mL_{[X,\alpha^r]} - \mL_{X_M} u^* \mL_{\alpha^r} + u^* \mL_{\alpha^r} \mL_{X}) F \\
& = (u^*\mL_{[X,\alpha^r]} - u^* \mL_{X} \mL_{\alpha^r} + u^* \mL_{\alpha^r} \mL_{X}) F \\
& = u^*(\mL_{[X,\alpha^r]} - \mL_{X} \mL_{\alpha^r} + \mL_{\alpha^r} \mL_{X}) F \\
& = u^*(\mL_{[X,\alpha^r]} - [\mL_{X} , \mL_{\alpha^r}]) F \\
& = 0. \qedhere
\end{align*}
\end{proof}


\section{The dgla module structure}
\label{sec: mod}
In this section we use the operation `$\bullet$' of Definition \ref{def: bull} and the identities proven in \S \ref{sec: ident} to construct on the cochain complex $\CmbG$ the structure of a dg-module over the dgla $\XmbG$.
In \S \ref{sec: morita} we prove that this module structure is Morita invariant in an appropriate sense, and so defines an object on the associated stack.

\subsection{The module structure}
The definitions of $\alpha \bullet F$, $X \bullet f$ and $X \bullet F$ in Definition \ref{def: bull} are linear in each variable and so determine linear maps $\GA \otimes \CmG \to \CM$, $\XmG \otimes \CM \to \CM$ and $\XmG \otimes \CmG \to \CmG$. 

\begin{thm}
\label{thm: mu}
The map 
\begin{align*}
\label{eqn: mu}
\mu:\XmbG \otimes \CmbG & \xrightarrow{} \CmbG \\
(\alpha , X) \otimes (f , F) & \mapsto (u^*\mL_{\alpha^r} F + \mL_{X_M} f , \mL_X F) \\
& = (\alpha \bullet F + X \bullet f , X \bullet F)
\end{align*}
makes $\CmbG$ into a differential graded module over the dgla $\XmbG$.
\end{thm}

\begin{proof}
The condition that $\mu$ is a chain map, see identities in \eqref{eqn:dglamod}, is equivalent to 
\[
\delta(x \bullet y) = \p x \bullet y + (-1)^{|x|} x \bullet \delta y
\]
for every homogeneous elements $x\in \XmbG$ and $y\in \CmbG$.
Indeed, this follows from Lemma \ref{lem: diff} once one takes into account the signs and uses the facts that $\p X = 0$, $\delta F = 0$ and $\alpha \bullet f = 0$ for $X \in \XmG$, $F \in \CmG$, $f \in \CG$ and $\alpha \in \GA$:
\begin{align*}
\delta ( X \bullet f ) & = \p X \bullet f + X \bullet \delta f = X \bullet \delta f & \text{ by Lemma \ref{lem: diff} \eqref{it: lem_diff_dXf}} \\
\delta ( X \bullet F ) & = \p X \bullet F + X \bullet \delta F & \text{as all terms have degree 2} \\
\delta ( \alpha \bullet f ) & = 0 = \p \alpha \bullet f - \alpha \bullet \delta f & \text{ by Lemma \ref{lem: diff} \eqref{it: lem_diff_adf}} \\
\delta ( \alpha \bullet F ) & = \p \alpha \bullet F - \alpha \bullet \delta F = \p \alpha \bullet F & \text{ by Lemma \ref{lem: diff} \eqref{it: lem_diff_daF}.}
\end{align*}

The condition that the action $\mu$ is compatible with the graded Lie brackets on $\XmbG$ is that 
\begin{equation}
\label{eqn: grLiebracket}
[x,x'] \bullet y = x \bullet (x' \bullet y) - (-1)^{|x||x'|} x' \bullet (x \bullet y)
\end{equation}
holds for $x$ a homogeneous element of $\XmbG$ and $y$ a homogeneous element of $\CmbG$. 
If $\alpha , \beta \in \GA$ then $\LL \alpha , \beta \RR(-), \alpha \bullet( \beta-)$ and $\beta\bullet(\alpha-)$ are operators of degree -2, thus 
\[
\LL \alpha , \beta \RR (y) = 0 = \alpha \bullet (\beta \bullet y) - (-1)^{|\alpha||\beta|} \beta \bullet (\alpha \bullet y)
\]
for any element $y \in \CmbG$. 
The remaining cases follow from Lemma \ref{lem: lie} \eqref{it: lem_lie_XYf}-\eqref{it: lem_lie_XaF} and the fact that multiplicative vector fields have degree zero.
\end{proof}

\subsection{Morita invariance and quasi-isomorphisms}
\label{sec: morita}
Recall from \S \ref{sec: dgalgebra} the categories $\dglmod$ (Definition \ref{def: dglmod}) and $\glmod$ (Definition \ref{def: glmod}).
The module structure of Theorem \ref{thm: mu} associates to a Lie groupoid $G$ an object $(\XmbG,\CmbG)$ in $\boldsymbol{\mathrm{dglmod}}$.
The setup and notation in the proof of the following result, as well as the statement itself, will be used a number of times in the sequel.

\begin{thm}
\label{thm: morita}
If $G$ and $H$ are Morita equivalent Lie groupoids then $(\XmbG,\CmbG)$ and $(\XmbH , C_m^{\bullet}(H))$ are quasi-isomorphic objects in the category $\boldsymbol{\mathrm{dglmod}}$.
\end{thm}

\begin{proof}
The proof is a combination of results showing the Morita invariance of differentiable cohomology \cite{Crainic03} and of the dgla $\XmbG$ \cite{OrtizW19, BerwickEvansL20}.
By a standard argument we can reduce the problem to the case where there exists a surjective Morita map $\phi : G \to H$.
Following \cite[\S 7.5]{OrtizW19} there exists a differential graded Lie subalgebra $\XmbG^{\phi} \subseteq \XmbG$ consisting of sections of the Lie algebroid $A$ of $G$ and multiplicative vector fields that are $\phi$-projectable.
Projecting elements of $\XmbG^{\phi}$ to the corresponding (uniquely defined) elements of $\XmbH$, and the inclusion of $\XmbG^{\phi}$ in $\XmbG$, define quasi-isomorphisms of dgla's
\begin{equation}
\label{eqn: qisom_XmbG}
\XmbG \xleftarrow{\mathrm{inc}} \XmbG^{\phi} \xrightarrow{\overline{\phi}} \XmbH .
\end{equation}
In loc.\ cit.\ the object $\XmbG^{\phi}$ is denoted $C_{\mathrm{mult}}^{\bullet} (\mathcal{V})^{\phi}$ and the statements are written in terms of crossed modules; see \S \ref{sec: 2lie} for a dictionary.
Following \cite[\S 1.2]{Crainic03}, pulling back functions by $\phi$ defines a quasi-isomorphism from $\CdG$ to $\CdH$, which we can truncate to give the quasi-isomorphism $\phi^*$ in the diagram
\begin{equation}
\label{eqn: qisom_CmbG}
\CmbG \xrightarrow{\mathrm{id}} \CmbG \xleftarrow{\phi^*} \CmbH .
\end{equation}
We can restrict the $\XmbG$-module $\CmbG$ to $\XmbG^{\phi}$ to define an object $(\XmbG^{\phi},\CmbG)$ in $\dglmod$.
We claim that combining \eqref{eqn: qisom_XmbG} and \eqref{eqn: qisom_CmbG} gives a diagram of quasi-isomorphisms in $\dglmod$
\begin{equation}
\label{eqn: qisom_dglmod}
( \XmbG , \CmbG ) \xleftarrow{(\mathrm{inc},\mathrm{id})} 
( \XmbG^{\phi} , \CmbG ) \xrightarrow{(\overline{\phi},\phi^*)}
( \XmbH , \CmbH ) 
\end{equation}
showing that $(\XmbG,\CmbG)$ and $(\XmbH,\CmbH)$ are quasi-isomorphic objects.
It is immediate that $(\mathrm{inc},\mathrm{id})$ is a morphism in $\dglmod$.
That $(\overline{\phi},\phi^*)$ is a morphism in $\dglmod$ follows from basic properties of vector fields and smooth maps: if $X \in \XmG$ is $\phi$-projectable to $\overline{\phi}(X) \in \XmH$ then for any function $F \in \CG$ it holds that 
\[
\phi^* \left(\mL_{\overline{\phi}(X)} F\right) = \mL_X (\phi^* F) .
\]
This, and the same identity applied to sections of the Lie algebroid of $G$, shows that $\phi^*$ defines a morphism of $\XmbG$ modules $\phi^* : \phi^* \CmbH \to \CmbG$.
By the results referenced above the morphisms $(\mathrm{inc},\mathrm{id})$ and $(\overline{\phi},\phi^*)$ are both quasi-isomorphisms.
\end{proof}

As a consequence of Theorem \ref{thm: morita}, one has the following.

\begin{cor}
\label{cor: morita_glmod}
In the setup and notation in the proof of Theorem \ref{thm: morita}.
The maps
\[
( \HXmbG , \HCmbG ) \xleftarrow{( H(\mathrm{inc}) , H(\mathrm{id}) )} 
( \HXmbG^{\phi} , \HCmbG ) 
\]
and
\[
( \HXmbG^{\phi} , \HCmbG ) \xrightarrow{( H(\overline{\phi}) , H(\phi^*) )}
( \HXmbH , \HCmbH ) 
\]
are isomorphisms in the category $\glmod$.
\end{cor}

\subsection{Quasi-isomorphisms of complexes of $\CMG$-modules}
\label{sec: CMG_qisom}
Recall from Definition \ref{def: CMG_mod} that $\XmbG$ is a complex of $\CMG$-modules, with action defined by $f \cdot \alpha = f\alpha$ and $f \cdot X = (t^*f)X$ for $f \in \CMG$, $\alpha \in \GA$ and $X \in \XmG$.

\begin{remi}
If $H \rightrightarrows N$ is a second Lie groupoid and $\phi : G \to H$ is a Morita map covering $\phi_0 : M \to N$ then $\phi^* : \CmbH \to \CmbG$ is a quasi-isomorphism and the map $\phi_0^* : \CNH \to \CMG$ is an isomorphism of algebras.
In particular, if $f \in \CMG$ then $f = \phi_0^* f'$ for a unique $f' \in \CNH$.
Via this isomorphism any complex of $\CNH$-modules can be considered as a complex of $\CMG$-modules.
\end{remi}

In the statement of the following two results we follow the setup and notation used in the proof of Theorem \ref{thm: morita}.
That is, $\phi : G \to H$ is a surjective Morita equivalence and $\XmbG^{\phi} \subseteq \XmbG$ is the subcomplex of $\phi$-projectable multiplicative vector fields and sections of $A$.

\begin{prop}
\label{prop: CMG_qisom}

Consider $\CmbH$ as a complex of $\CMG$-modules via the algebra isomorphism $\phi_0^* : \CNH \to \CMG$. Then:
\begin{enumerate}
\item \label{it: CMG_subcomplex} 
$\XmbG^{\phi} \subseteq \XmbG$ is a subcomplex of $\CMG$-modules.

\item \label{it: CMG_morphisms}
The maps
\[
\XmbG \xleftarrow{\mathrm{inc}} \XmbG^{\phi} \xrightarrow{\overline{\phi}} \XmbH \]
are quasi-isomorphisms of complexes of $\CMG$-modules.
\end{enumerate}
\end{prop}

\begin{proof}
\eqref{it: CMG_subcomplex}.
Suppose that $f = \phi_0^*f'$ is an element of $\CMG$ and $X \in \XmG$ is projectable to a multiplicative vector field $X'$ on $H$.
Then $(t^*f)X$ is projectable to $(t^*f')X'$ because if $k \in C^{\infty}(H)$ then 
\begin{align*}
\mL_{ (t^*f)X } (\phi^*k) & = (t^*f) \mL_{X} (\phi^*k) \\
& = (t^*\phi_0^*f') \mL_{ X} (\phi^*k) \\
& = (\phi^*t^*f') \mL_{ X } (\phi^*k) \\
& = (\phi^*t^*f') \phi^* \left(\mL_{ X' } k \right) \\
& = \phi^*  \left(  (t^*f') \left(\mL_{ X' } k \right) \right) \\
& = \phi^*  \left(\mL_{ (t^*f') X' } k \right)
\end{align*}
where in the third equality we have used the fact that $\phi_0 t = t \phi$ because $\phi : G \to H$ is a morphism of Lie groupoids.
Similarly, if $\alpha \in \GA$ is projectable to $\alpha' \in \GB$ then $f\alpha$ is projectable to $f' \alpha'$ because
\[
(f\alpha)(x) = \left( (\phi_0^*f')\alpha\right)(x) = f'(\phi_0(x)) \alpha(x)
\]
so that 
\[
\mathrm{Lie}(\phi)  \left( (f\alpha)(x) \right) = (f'\alpha') (\phi_0(x))
\]
where $\mathrm{Lie}(\phi) : A \to B$ is the morphism of Lie algebroids associated to $\phi$.
It follows that $\XmbG^{\phi}$ is closed under the action of $\CMG$ and so is a subcomplex of $\CMG$-modules.

\eqref{it: CMG_morphisms}.
By statement \eqref{it: CMG_subcomplex}, $\XmbG^{\phi}$ is a subcomplex of $\CMG$-modules and so the inclusion $\XmbG^{\phi} \hookrightarrow \XmbG$ is automatically a morphism of complexes of $\CMG$-modules.
Both this inclusion and the map $\overline{\phi} : \XmbG^{\phi} \to \XmbH$ are quasi-isomorphisms of complexes of vector spaces by the results referenced in the proof of Theorem \ref{thm: morita}.
It remains to prove that each component of $\overline{\phi}$ is a morphism of $\CMG$-modules.

With respect to the isomorphism $\phi_0^* : \CNH \to \CMG$ an element $f \in \CMG$ acts on $\XmbH$ via $f' \in \CNH$, where $f'$ is the unique element of $\CNH$ satisfying $\phi_0^*f' = f$.
The statement then follows from the computations in the proof of part \eqref{it: CMG_subcomplex}:
\begin{multicols}{2}
\begin{align*}
\overline{\phi} ( f \cdot X ) & = \overline{\phi} ( (t^*f)X ) \\
& = (t^*f') \overline{\phi} (X) \\
&=f'\cdot \overline{\phi} (X)\\
& = f \cdot \overline{\phi} (X), \\
\end{align*}

\begin{align*}
\overline{\phi} ( f \cdot \alpha ) & = \overline{\phi} ( f\alpha ) \\
& = f' \overline{\phi} (\alpha) \\
& = f \cdot \overline{\phi} (\alpha) .
\qedhere
\end{align*}
\end{multicols}

\end{proof}

\begin{cor}
\label{cor: CMG_qisom_grmod}
Consider $\CmbH$ as a complex of $\CMG$-modules via the algebra isomorphism $\phi_0^* : \CNH \to \CMG$.
The maps
\[
\HXmbG \xleftarrow{ H(\mathrm{inc}) } H^{\bullet} \left((\XmbG)^{\phi}\right) \xrightarrow{ H(\overline{\phi}) } \HXmbH
\]
are isomorphisms of graded $\CMG$-modules.
\end{cor}

Combining all the previous results, one gets the following.

\begin{prop}
\label{prop: CMG_qisom_general}
If $G$ and $H$ are Morita equivalent Lie groupoids then there is an algebra isomorphism $\CNH \to \CMG$ with respect to which $\XmbG$ and $\XmbH$ are quasi-isomorphic as complexes of $\CMG$-modules.
\end{prop}

\subsection{Objects on stacks}
The results in \S \ref{sec: mod} and \S \ref{sec: morita} can be interpreted in terms of differentiable stacks as follows.
Recall from Proposition \ref{prop: fun} that if $\mathcal{X}$ is a differentiable stack then the algebra $\CX$ is intrinsically defined, and is isomorphic to $\CM^G$ if $G \rightrightarrows M$ is the Lie groupoid associated to an atlas $M \to \mathcal{X}$.

\begin{thm}
\label{thm: dglmod_stacks}
Let $\mathcal{X}$ be a differentiable stack and $G \rightrightarrows M$ the Lie groupoid associated to an atlas $M \to \mathcal{X}$. Then, the object 
\[
(\XmbX,\CmbX) := (\XmbG,\CmbG)
\]
is well defined in the category $\boldsymbol{\mathrm{dglmod}}$ and $\XmbX$ is a complex of $\CX$-modules.
Up to quasi-isomorphism these objects are independent of the choice of atlas and unchanged if $\mathcal{X}$ is replaced by an equivalent stack $\mathcal{Y}$.
\end{thm}

\begin{proof}
This follows from Theorem \ref{thm: morita} and Proposition \ref{prop: CMG_qisom_general}, and the fact that if $M \to \mathcal{X}$ and $N \to \mathcal{X}$ are atlases of $\mathcal{X}$, and $M' \to \mathcal{Y}$ is an atlas of an equivalent stack $\mathcal{Y}$, then the three associated Lie groupoids are Morita equivalent.
\end{proof}


\section{Cohomology and derivations}
\label{sec: co}

In this section we show that the graded Lie algebra $\HXmbG$ carries a natural graded Lie-Rinehart algebra structure (see Definition \ref{def: grLR}) over the algebra $\CM^G = H^{0}(\CmbG)$. As a consequence, associated to each differentiable stack $\mathcal{X}$ is a graded Lie-Rinehart algebra over the algebra $\CX$, independent up to isomorphism of the choice of atlas involved in the construction.

\subsection{Some module structures}
Recall from  Theorem \ref{thm: mu} that there is a dgla $\XmbG$-module structure on $\CmbG$ given by:
\begin{align*}
\mu : \XmbG \otimes \CmbG & \to \CmbG \\
(\alpha , X) \otimes (f , F) & \mapsto (u^*\mL_{\alpha^r} F + \mL_{X_M} f , \mL_X F) .
\end{align*}
The corresponding map on cohomology $H(\mu)$ defines a graded $\HXmbG$-module structure on $\HCmbG$:
\begin{align}
\label{eqn: Hmu}
H(\mu) : H^{\bullet} (\XmbG) \otimes H^{\bullet} (\CmbG) & \to H^{\bullet} (\CmbG) \\
\nonumber
( [\alpha] , [X] ) \otimes ( [f] , [F] ) & \mapsto ( [u^* \mL_{\alpha^r} F + \mL_{X_M} f] , [\mL_X F] ) 
\end{align}
where the K\"{u}nneth theorem has been applied on the left side.
The map $H(\mu)$ corresponds to a morphism of dgla's:
\begin{align}
\label{eqn: XmbG_End}
\bmL : H^{\bullet} (\XmbG) & \to \mathrm{End}^{\bullet} (\HCmbG \\
\nonumber
x & \mapsto \left(y \mapsto H(\mu) (x \otimes y) \right).
\end{align}

We also have a left graded $\HCmbG$-module structure on $\HXmbG$ introduced in Definition \ref{def: cdot}:
\begin{align}
\label{eqn: cdot2}
\rho : \HCmbG \otimes \HXmbG & \to \HXmbG \\
\nonumber
\left( [f],[F] \right) \otimes \left( [\alpha],[X] \right) & \mapsto \left( [f \alpha] , [F \alpha^r] + [(t^*f) X] \right) .
\end{align}

\begin{remi} If $G$ is the unit groupoid $M \rightrightarrows M$ then $\HXmbG$ is equal to $\XM$, $\HCmbG$ coincides with $\CM$, $\bmL(X) (f) = \mL_{X} f$ and $\rho(f \otimes X) = fX$. In other words, the map  $\bmL$ \eqref{eqn: XmbG_End} is just the Lie derivative and $\rho$ \eqref{eqn: cdot2} is the standard $\CM$-module structure on $\XM$.
In this case $\bmL$ is an isomorphism onto the Lie algebra $\mathrm{Der}(\CM)$
and the Lie bracket on $\XM$ satisfies the Leibniz identity.
\end{remi}

The main results in the present section \S \ref{sec: co} are analogues for more general Lie groupoids and for differentiable stacks.

We start by establishing some analogues of the Leibniz rule and proving that under certain conditions the image of the morphism of gla's $\bmL$ is contained in the gla $\mathrm{Der}(\HCmbG)$ of graded derivations of the graded algebra $\HCmbG$. Recall from \S \ref{sec: CmbG} that if $G$ is proper then $\HCmbG = H^{0}(\CmbG)$, in which case the results in this section can be restated with $\HCmbG$ in place of $H^{0}(\CmbG)$ (or with $\HCmbX$ in place or $\CX$ for $\mathcal{X}$ the stack corresponding to $G$).
For the non-proper case see \S \ref{ex: non-prop}.

\subsection{The graded Lie-Rinehart structure}
The dgla $\XmbG$-module structure on $\CmbG$ descends to cohomology to define a gla $\HXmbG$-module structure on $\HCmbG$. The commutative algebra $\CM^G$ is the degree zero graded submodule of $H^{\bullet}(\CmbG)$, so that the restriction map
\begin{align*}
H^{\bullet} (\XmbG) \otimes \CM^G & \to \CM^G \\
( [\alpha] , [X] ) \otimes f & \mapsto \alpha \bullet f + X \bullet f = \mL_{X_M} f
\end{align*}
makes $\CM^G$ into a module over the graded Lie algebra $\HXmbG$.
This module structure determines a morphism of graded Lie algebras
\begin{equation}\label{eq:anchor}\bmL : \HXmbG \to \mathrm{Der} \left( \CM^G \right), \quad ([\alpha],[X]) \mapsto \mathcal{L}_{X_M}\end{equation}
where we have restricted the codomain to the Lie subalgebra $ \mathrm{Der} \left( \CM^G \right)$ of $\mathrm{End} \left( \CM^G \right)$ and $\mathcal{L}_{X_M}$ is a derivation because it is a Lie derivative.
Recall from Definition \ref{def: CMG_mod} the left $\CM^G$-module structure on $\HXmbG$ defined by $f \cdot [\alpha] = [f\alpha]$ and $f \cdot [X] = [(t^*f) X]$.

\begin{prop}
\label{prop: grLR}
The $\CM^G$-module structure on $\HXmbG$ and the morphism of graded Lie algebras $\bmL$ make $\HXmbG$ into a graded Lie-Rinehart algebra over $\CM^G$.
That is:
\begin{enumerate}
\item \label{it: prop_grLR_bmL_morph}
The map $\bmL$ is both a morphism of dgla's and $\CM^G$-modules,

\item \label{it: prop_grLR_grLeib}
The graded Lie bracket on $\HXmbG$ satisfies the graded Leibniz identity.
\end{enumerate}
\end{prop}

\begin{proof}
We start by proving \emph{(a)}. From equation \ref{eq:anchor} follows that $\bmL$ is a morphism of dgla's.
Suppose that $f \in \CM^G$ and $X \in \XmG$ is a multiplicative vector field on $G$.
Then $(t^*f)X$ is a multiplicative vector field with $((t^*f)X)_M = fX_M$ and as derivations of $\CM^G$
\[
\mL_{ ((t^*f)X)_M } = \mL_{ fX_M } = f \mL_{X_M} .
\]
This shows that $\bmL$ is a morphism of $\CM^G$-modules.

Regarding \emph{(b)}, suppose that $f \in \CM^G$, $\alpha \in \Gamma(A)$ is a cocycle and $X,Y \in \XmG$ are multiplicative vector fields.
Then
\begin{align*}
[ X , (t^*f)Y ] & = \mL_{X}(t^*f) Y + f[X,Y] \\
& = t^* (\mL_{X_M}f) Y + f[X,Y]
\end{align*}
and
\begin{align*}
\LL X , f\alpha \RR & = [ X , (t^*f)\alpha^r ]|_{M} \\
& = \left( \mL_{X}(t^*f) \alpha^r + (t^*f) [X,\alpha^r] \right)|_{M} \\
& = \left( t^*(\mL_{X_M}f) \alpha^r + (t^*f) [X,\alpha^r] \right)|_{M} \\
& = (\mL_{X_M}f) \alpha + f \LL X,\alpha \RR .
\end{align*}
Passing to cohomology then shows that the graded Lie bracket on $\HXmbG$ satisfies the graded Leibniz identity.
\end{proof}

\vspace{.2cm}
\noindent \textbf{Morita invariance and stacks.}
Using the Morita invariance result of Theorem \ref{thm: morita} we can prove a similar statement for the Lie-Rinehart algebra of Proposition \ref{prop: grLR}. For that, we follow the setup and notation in the proof of Theorem \ref{thm: morita}, Corollary \ref{cor: morita_glmod} and Corollary \ref{cor: CMG_qisom_grmod}.
That is, $\phi : G \to H$ is a surjective Morita equivalence and $\XmbG^{\phi} \subseteq \XmbG$ is the subcomplex of $\phi$-projectable multiplicative vector fields and sections of $A$. Identifying $\CNH$ with $\CMG$ via the algebra isomorphism $\phi_0^* : \CNH \to \CMG$, one has the following result.

\begin{prop}
\label{prop: LR_isom}
The map
\[
 \HXmbG \xrightarrow{ H(\overline{\phi}) \circ H(\mathrm{inc})^{-1} }
\HXmbH 
\]
is an isomorphism of graded Lie-Rinehart algebras over $\CMG$.
\end{prop}

\begin{proof}
By Corollary \ref{cor: morita_glmod} and Corollary \ref{cor: CMG_qisom_grmod} the maps $H(\mathrm{inc})$ and $H(\overline{\phi})$ are isomorphisms of graded Lie algebras and of graded $\CMG$-modules.
After identifying $\CNH$ with $\CMG$ the isomorphisms in the category $\glmod$ given in Corollary \ref{cor: morita_glmod} show that $H(\mathrm{inc})$ and $H(\overline{\phi})$ are compatible with the actions of $\HXmbG$ and $\HXmbH$ on $\CMG$.
\end{proof}

Recall from Proposition \ref{prop: fun} that if $\mathcal{X}$ is a differentiable stack then the algebra $\CX$ is intrinsically defined, and is isomorphic to $\CM^G$ if $G \rightrightarrows M$ is the Lie groupoid associated to an atlas $M \to \mathcal{X}$. In this setting, the following holds.

\begin{thm}
\label{thm: co_stacks}
There is an associated a graded Lie-Rinehart algebra 
\[
\mathfrak{X}^{\bullet}(\mathcal{X}) := \HXmbG
\]
over the algebra $\CX$, whose isomorphism class is independent of the choice of atlas, making it a well-defined object in the category of $\mathbf{gLRalg}$. Also, if $\mathcal{X}$ is equivalent to $\mathcal{Y}$, then $\mathfrak{X}^{\bullet}(\mathcal{X})\cong \mathfrak{X}^{\bullet}(\mathcal{Y})$ as graded Lie-Rinehart algebras over $\CX$. 
\end{thm}

\begin{proof}
This follows from Proposition \ref{prop: LR_isom} and the fact that if $M \to \mathcal{X}$ and $N \to \mathcal{X}$ are atlases of $\mathcal{X}$, and $M' \to \mathcal{Y}$ is an atlas of an equivalent stack $\mathcal{Y}$, then the three associated Lie groupoids are Morita equivalent.
\end{proof}

\begin{remi}
A natural question is whether or not the morphism 
\begin{equation} \label{eq:isomorphism}
    \mathfrak{X}^{\bullet}(\mathcal{X}) \to \mathrm{Der}(\CX)
\end{equation}
is an isomorphism as in the case of a smooth manifold.
In \S \ref{sec: ex} we study a number of examples and show that in general the answer to this question is no, an obvious obstruction being the possible non-vanishing of $H^{-1}(\XmbG)$ where $G$ is a Lie groupoid with classifying stack equivalent to $\mathcal{X}$.
\end{remi}


\section{The infinitesimal picture}
\label{sec: inf}
Let $G$ be a Lie groupoid over $M$ with Lie algebroid $A$. Given some cochain complex associated to $G$, a Van-Est type map relates such a complex to a complex defined purely in terms of $A$. This is the case of the Van-Est maps $\oVE$ (Definition \ref{def: VE}) and $\VE$ (Definition \ref{def: oVE}). The following results says that these Van-Est maps are also compatible with the module structures $\overline{\mu}$ (Definition \ref{def: obull}) and $\mu$ (Theorem \ref{thm: mu}).

\begin{thm}
\label{thm: VE}
The following statements hold:
\begin{enumerate}
\item \label{it: thm_VE_dgla}
The Van-Est map $\VE : \XmbG \to \XmbA$ is a morphism of dgla's.
\item \label{it: thm_VE_mod} The following diagram is a commutative diagram of morphisms of cochain complexes:
\begin{equation}
\label{diag: VE_mod}
\xymatrix{
\XmbG \otimes \CmbG \ar[d]_{\VE \otimes \oVE} \ar[r]^-{\mu} & \CmbG \ar[d]^{\oVE} \\
\XmbA \otimes \CmbA \ar[r]_-{\overline{\mu}} & \CmbA
}
\end{equation}
\item \label{it: thm_VE_isom} The  vertical arrows \eqref{diag: VE_mod} are isomorphisms whenever $G$ is source simply connected.
\end{enumerate}
\end{thm}

\begin{proof}
Note that \emph{(a)} follows directly from Remark \ref{rem: adrep}. Regarding \emph{(b)}, the maps $\oVE$, $\VE$, $\overline{\mu}$ and $\mu$ are chain maps by statement \eqref{it: thm_VE_dgla}, Definition \ref{def: oVE}, Definition \ref{def: obull} and Theorem \ref{thm: mu} respectively; it follows that $\VE \otimes \oVE$ is a chain map also.
It remains to show that the diagram \eqref{diag: VE_mod} commutes, which is equivalent to the statement that
\begin{equation*}
\oVE \, (x \bullet y) = \VE(x) \, \obullet \, \oVE (y).
\end{equation*}
for all simple tensors $x \otimes y \in \XmbG \otimes \CmbG$ with $x$ and $y$ homogeneous.
We deal with the possible cases in turn, using the Definitions \ref{def: obull}, \ref{def: oVE}, \ref{def: VE} and \ref{def: bull}.

If $\alpha \in \GA$ and $f \in \CM$ then
\begin{align*}
\oVE (\alpha \bullet f) - \VE (\alpha) \obullet \oVE (f) 
& = \oVE(0) - \alpha \obullet f \\
& = 0. 
\end{align*}

If $\alpha \in \GA$ and $F \in \CmG$ then
\begin{align*}
\oVE (\alpha \bullet F) - \VE(\alpha) \obullet \oVE(F) 
& = \oVE(u^*\mL_{\alpha^r}F) - \alpha \obullet \omega_F \\
& = u^*\mL_{\alpha^r}F - i_\alpha \omega_F \\
& = u^*\mL_{\alpha^r}F - u^*\mL_{\alpha^r} F \\
& = 0.
\end{align*}

If $X \in \XmG$ and $f \in \CM$ then
\begin{align*}
\oVE (X \bullet f) - \VE(X) \obullet \oVE(f) 
& = \oVE (\mL_{X_M} f) - D_X \obullet f \\
& = \mL_{X_M} f - \mL_{X_M} f \\
& = 0.
\end{align*}

If $X \in \XmG$ and $f \in \CmG$ then
\begin{align}
\nonumber
\oVE (X \bullet F) - \VE(X) \obullet \oVE(F) 
& = \oVE (\mL_X F) - D_X \obullet \omega_F \\
\label{eqn: VE_XF}
& = \omega_{\mL_X F} - (\mL_{X} \circ \omega_{F} - \omega_F \circ D_X ).
\end{align}
The right hand side of \eqref{eqn: VE_XF} is the section of $A^*$ given by the map
\begin{align*}
\alpha & \mapsto \omega_{\mL_X F} (\alpha)  - \mL_{X_M} (\omega_F (\alpha)) + \omega_F (D_X (\alpha)) \\
& = u^* (\mL_{\alpha^r} \mL_X F) - \mL_{X_M} (u^* \mL_{\alpha^r} F) + u^* (\mL_{[X,\alpha^r]} F) \\
& = u^* (\mL_{\alpha^r} \mL_X F) - u^* (\mL_X \mL_{\alpha^r} F) + u^* (\mL_{[X,\alpha^r]} F) \\
& = u^* ( \mL_{\alpha^r} \mL_X F - \mL_X \mL_{\alpha^r} F +  \mL_{[X,\alpha^r]} F ) \\
& = u^* ( \mL_{[\alpha^r,X]} - \mL_{[\alpha^r,X]} F ) \\
& = 0
\end{align*}
where in the second equality we have used the fact that that $X_M$ is $u$-related to $X$, so that $u^*\mL_X = \mL_{X_M} u^*$.

Finally, \emph{(c)} follows from the fact that the vertical arrows in \eqref{diag: VE_mod} are isomorphisms whenever $G$ is source simply connected follows from the properties of the Van-Est maps discussed in \S \ref{sec: van_est}.
\end{proof}


\section{Examples}
\label{sec: ex}

In this section we give examples and conditions under which the map \eqref{eq:isomorphism} is an isomorphism. 

\begin{ex} 

Let $G=(N\times_MN\rightrightarrows N)$ be the submersion groupoid defined by a surjective submersion $\pi:N\to M$. The corresponding stack $\mathcal{X}$ is equivalent to the manifold $M$. Hence the graded Lie-Rinehart structure on $\mathfrak{\mathcal{X}}$ is equivalent to the standard Lie-Rinehart algebra structure on the base $M$

 $$\mathcal{L}:\mathfrak{X}(M)\to \mathrm{Der}(C^{\infty}(M)).$$
 
\noindent In this case, the map \eqref{eq:isomorphism} is an isomorphism. 
 
\end{ex}

Let $G\rightrightarrows M$ be a proper Lie groupoid with Lie algebroid $A$. There are well-defined sets $\Gamma(\mathrm{i}):= \ker(\rho: \Gamma(A)\to \mathfrak{X}(M))$ and $\Gamma(\nu):=\mathfrak{X}(M)/\mathrm{im}(\rho)$. By Theorem 6.1 in \cite[\S6]{CrainicMS20} one has
\[
H^{\bullet}(\mathfrak{X}_{m}^{\bullet}(G))\simeq \Gamma(\mathrm{i})^{inv}\oplus \Gamma(\nu)^{inv},
\]
\noindent where $\Gamma(\mathrm{i})^{inv}$ is the set of invariant sections by the adjoint action of $G$, and  $[V] \in \Gamma(\nu)$ is invariant if there exists $X\in \mathfrak{X}(G)$ which is both $t$-projectable and $s$-projectable to $V$, see \S 4 in \cite{CrainicMS20} for more details.  On the other hand, by Proposition 1 in \cite[\S 2.1]{Crainic03} the cohomology of $G$ is concentrated in degree zero, hence
\[
H^{\bullet}(C_m^{\bullet}(G))=C^{\infty}(M)^{G}\oplus 0.
\]
We will apply this to the next example.

\begin{ex}
Let $\pi:T\to M$ be a bundle of tori. This defines a proper Lie groupoid $\mathrm{T}\rightrightarrows M$ with both source and target maps given by $\pi$. By the previous observations we get that the morphism \eqref{eq:isomorphism} is given by
\[\Gamma(A)\oplus \mathfrak{X}(M)\to \mathrm{Der}(C^{\infty}(M)), \]
where $A=\ker(T\pi)$ is a Lie algebroid of $\mathrm{T}\rightrightarrows M$. In particular, if $A\neq 0$, then the morphism \eqref{eq:isomorphism} is not an isomorphism.
\end{ex}

If $G$ is an \'etale Lie groupoid then $A = 0$ and the map \eqref{eq:isomorphism} is the map
\[
\mathcal{L} : \mathfrak{X}(M)^{G} \to \mathrm{Der}(C^{\infty}(M)^{G})
\]
given by the action of invariant vector fields on invariant functions.
The following example shows that even for proper and \'etale groupoids, i.e.\ those presenting orbifolds, the map \eqref{eq:isomorphism} is in general not an isomorphism.

\begin{ex}
Let $\mathbb{Z}_2$ act on $\mathbb{R}$ by reflection around zero. Then the orbit space $\mathbb{R}/\mathbb{Z}_2$ is a global orbifold presented by the action groupoid $\mathbb{Z}_2 \ltimes \mathbb{R} \rightrightarrows \mathbb{R}$. The space of invariant functions $C^{\infty}(\mathbb{R})^{\mathbb{Z}_2}$ consists of all even functions, and the space of invariant vector fields $\mathfrak{X}(\mathbb{R})^{\mathbb{Z}_2}$ consists of all vector fields $X = \phi(x) \frac{\partial}{\partial x}$ where $\phi(x)$ is an odd function.
It follows from Hadamard's Lemma that the map $D := \frac{1}{x} \frac{\partial}{\partial x}$ is a well-defined algebra derivation of $C^{\infty}(\mathbb{R})^{\mathbb{Z}_{2}}$.
If $X$ is an invariant vector field then $\mathcal{L}_{X}(x^2)$ vanishes at zero, whereas $D(x^2) = 2$, and therefore the map \eqref{eq:isomorphism} $\mathfrak{X}(\mathbb{R})^{\mathbb{Z}_2} \to \mathrm{Der}(C^{\infty}(\mathbb{R})^{\mathbb{Z}_2})$ is not surjective.
\end{ex}

The following example was motivated by the work of P. Molino in \cite[\S 2.2]{Molino88}. In the example below we show that the map \eqref{eq:isomorphism} vanishes identically.

\begin{ex}[Molino's example]\label{ex: non-prop}
Let $\mathbb{T}^2=S^1\times S^1$ be the torus with local coordinates $(\theta_0,\theta_1)$. Consider the foliation $\mathcal{F}$ on $\mathbb{T}^2$ determined by the tangent distribution $D=\mathrm{span}\{\lambda_0\frac{\partial}{\partial \theta_0}+\lambda_1\frac{\partial}{\partial \theta_1}\mid \lambda_1/\lambda_0 \in \mathbb{R}\setminus \mathbb{Q}\}$. Let $G\rightrightarrows \mathbb{T}^2$ be the foliation groupoid whose canonical foliation is $\mathcal{F}$. On the one hand, it is clear that $G$ is not a proper groupoid and $C^{\infty}(\mathbb{T}^2)^{G}$ consists of constant functions, which implies that $\mathrm{Der}(C^{\infty}(\mathbb{T}^2)^{G})=0.$ On the other hand, by Proposition 3.3 in \cite[\S 3]{CrainicMS20} we have that $H^{\bullet}(\mathfrak{X}_m^{\bullet}(G))\simeq 0\oplus \mathrm{span}\{-\lambda_1\frac{\partial}{\partial \theta_0}+\lambda_0\frac{\partial}{\partial \theta_1}\}$. Thus, in this case the map \eqref{eq:isomorphism} is the zero morphism, which is not an isomorphism.
\end{ex}

\begin{ex}
Let $G$ be a Lie group considered as a Lie groupoid over a point, and $BG$ the classifying stack of $G$.
Then $\HXmbG = \mathfrak{g}^{G} \oplus H^{1}_{d}(G,\mathfrak{g})$, $C^{\infty}(BG) = \mathbb{R}$, and $\mathrm{Der}(\mathbb{R}) = 0$.
It follows that the map \eqref{eq:isomorphism} is an isomorphism iff $\mathfrak{g}^{G} = 0$ and $H^{1}_{d}(G,\mathfrak{g}) = 0$.
The latter condition holds if $G$ is compact or semi-simple with finite fundamental group, but fails in general, e.g.\ for $G = GL_{2}(\mathbb{C})$.
Note that if $G$ is connected then $\mathfrak{g}^{G} = Z(\mathfrak{g})$, and if $G$ is moreover simply connected then $H^{1}_{d}(G,\mathfrak{g}) = H^{1}(\mathfrak{g},\mathfrak{g})$.
\end{ex}


\bibliography{derstack}{}

\providecommand{\bysame}{\leavevmode\hbox to3em{\hrulefill}\thinspace}
\providecommand{\MR}{\relax\ifhmode\unskip\space\fi MR }
\providecommand{\MRhref}[2]{%
  \href{http://www.ams.org/mathscinet-getitem?mr=#1}{#2}
}
\providecommand{\href}[2]{#2}
\begin{thebibliography}{BCLGX20}

\bibitem[AAC11]{AbadC11}
Camilo Arias~Abad and Marius Crainic, \emph{The {W}eil algebra and the {V}an
  {E}st isomorphism}, Ann. Inst. Fourier (Grenoble) \textbf{61} (2011), no.~3,
  927--970. \MR{2918722}

\bibitem[AAC13]{AbadC13}
\bysame, \emph{Representations up to homotopy and {B}ott's spectral sequence
  for {L}ie groupoids}, Adv. Math. \textbf{248} (2013), 416--452. \MR{3107517}

\bibitem[AAS11]{AbadS11}
Camilo Arias~Abad and Florian Sch\"{a}tz, \emph{Deformations of {L}ie brackets
  and representations up to homotopy}, Indag. Math. (N.S.) \textbf{22} (2011),
  no.~1-2, 27--54. \MR{2853613}

\bibitem[AC12]{AbadC12}
Camilo~Arias Abad and Marius Crainic, \emph{Representations up to homotopy of
  {L}ie algebroids}, J. Reine Angew. Math. \textbf{663} (2012), 91--126.
  \MR{2889707}

\bibitem[BC04]{BaezC04}
John~C. Baez and Alissa~S. Crans, \emph{Higher-dimensional algebra. {VI}. {L}ie
  2-algebras}, Theory Appl. Categ. \textbf{12} (2004), 492--538. \MR{2068522}

\bibitem[BCLGX20]{BonechiCLGX20}
Francesco Bonechi, Nicola Ciccoli, Camille Laurent-Gengoux, and Ping Xu,
  \emph{Shifted poisson structures on differentiable stacks}, 2020.

\bibitem[BEL20]{BerwickEvansL20}
Daniel Berwick-Evans and Eugene Lerman, \emph{Lie 2-algebras of vector fields},
  Pacific J. Math. \textbf{309} (2020), no.~1, 1--34. \MR{4202003}

\bibitem[BX03]{BehrendX03}
Kai Behrend and Ping Xu, \emph{{$S^1$}-bundles and gerbes over differentiable
  stacks}, C. R. Math. Acad. Sci. Paris \textbf{336} (2003), no.~2, 163--168.
  \MR{1969572}

\bibitem[CM08]{CrainicM08}
Marius Crainic and Ieke Moerdijk, \emph{Deformations of {L}ie brackets:
  cohomological aspects}, J. Eur. Math. Soc. (JEMS) \textbf{10} (2008), no.~4,
  1037--1059. \MR{2443928}

\bibitem[CMS20]{CrainicMS20}
Marius Crainic, Jo\~{a}o~Nuno Mestre, and Ivan Struchiner, \emph{Deformations
  of {L}ie groupoids}, Int. Math. Res. Not. IMRN (2020), no.~21, 7662--7746.
  \MR{4176835}

\bibitem[Cra03]{Crainic03}
Marius Crainic, \emph{Differentiable and algebroid cohomology, van {E}st
  isomorphisms, and characteristic classes}, Comment. Math. Helv. \textbf{78}
  (2003), no.~4, 681--721. \MR{2016690}

\bibitem[Dae14]{Daenzer14}
Calder Daenzer, \emph{Geometric {T}-dualization}, String-{M}ath 2013, Proc.
  Sympos. Pure Math., vol.~88, Amer. Math. Soc., Providence, RI, 2014,
  pp.~243--258. \MR{3330292}

\bibitem[Gri21]{Grivaux20}
Julien Grivaux, \emph{Derived intersections and free dg-{Lie} algebroids},
  Publ. Res. Inst. Math. Sci. \textbf{57} (2021), no.~3-4, 1049--1087
  (English).

\bibitem[Hep09]{Hepworth09}
Richard Hepworth, \emph{Vector fields and flows on differentiable stacks},
  Theory Appl. Categ. \textbf{22} (2009), 542--587. \MR{2591949}

\bibitem[KP23]{BjarneP2023}
Bjarne Kosmeijer and Hessel Posthuma, \emph{Lie groupoid deformations and
  convolution algebras}, J. Geom. Phys. \textbf{194} (2023), Paper No. 105012,
  23. \MR{4651862}

\bibitem[KV94]{KapranovV91}
M.~M. Kapranov and V.~A. Voevodsky, \emph{{$2$}-categories and {Z}amolodchikov
  tetrahedra equations}, Algebraic groups and their generalizations: quantum
  and infinite-dimensional methods ({U}niversity {P}ark, {PA}, 1991), Proc.
  Sympos. Pure Math., vol.~56, Amer. Math. Soc., Providence, RI, 1994,
  pp.~177--259. \MR{1278735}

\bibitem[Ler10]{Lerman10}
Eugene Lerman, \emph{Orbifolds as stacks?}, Enseign. Math. (2) \textbf{56}
  (2010), no.~3-4, 315--363. \MR{2778793}

\bibitem[Mac05]{MackenzieBook05}
Kirill C.~H. Mackenzie, \emph{General theory of {L}ie groupoids and {L}ie
  algebroids}, London Mathematical Society Lecture Note Series, vol. 213,
  Cambridge University Press, Cambridge, 2005. \MR{2157566}

\bibitem[Man04]{Manetti04}
Marco Manetti, \emph{Lectures on deformations of complex manifolds
  (deformations from differential graded viewpoint)}, Rend. Mat. Appl. (7)
  \textbf{24} (2004), no.~1, 1--183. \MR{2130146}

\bibitem[MM02]{MoerdijkM02}
Ieke Moerdijk and Janez Mr\v{c}un, \emph{On integrability of infinitesimal
  actions}, Amer. J. Math. \textbf{124} (2002), no.~3, 567--593. \MR{1902889}

\bibitem[MM03]{MoerdijkMBook03}
I.~Moerdijk and J.~Mr\v{c}un, \emph{Introduction to foliations and {L}ie
  groupoids}, Cambridge Studies in Advanced Mathematics, vol.~91, Cambridge
  University Press, Cambridge, 2003. \MR{2012261}

\bibitem[MM05]{MoerdijkMBook05}
\bysame, \emph{Lie groupoids, sheaves and cohomology}, Poisson geometry,
  deformation quantisation and group representations, London Math. Soc. Lecture
  Note Ser., vol. 323, Cambridge Univ. Press, Cambridge, 2005, pp.~145--272.
  \MR{2166453}

\bibitem[Mol88]{Molino88}
Pierre Molino, \emph{Riemannian foliations. {With} appendices by {G}. {Cairns},
  {Y}. {Carri{\`e}re}, {E}. {Ghys}, {E}. {Salem}, {V}. {Sergiescu}}, transl.
  from the {French} by {Grant} {Cairns} ed., Prog. Math., vol.~73, Boston, MA
  etc.: Birkh{\"a}user Verlag, 1988 (English).

\bibitem[MX98]{MackenzieX98}
Kirill C.~H. Mackenzie and Ping Xu, \emph{Classical lifting processes and
  multiplicative vector fields}, Quart. J. Math. Oxford Ser. (2) \textbf{49}
  (1998), no.~193, 59--85. \MR{1617335}

\bibitem[Nui19]{Nuiten19}
Joost Nuiten, \emph{Homotopical algebra for {L}ie algebroids}, Appl. Categ.
  Structures \textbf{27} (2019), no.~5, 493--534. \MR{4003715}

\bibitem[OW19]{OrtizW19}
C.~Ortiz and J.~Waldron, \emph{On the {L}ie 2-algebra of sections of an
  {$\mathcal{LA}$}-groupoid}, J. Geom. Phys. \textbf{145} (2019), 103474, 34.
  \MR{3985597}

\bibitem[Rin63]{Rinehart63}
George~S. Rinehart, \emph{Differential forms on general commutative algebras},
  Trans. Amer. Math. Soc. \textbf{108} (1963), 195--222. \MR{154906}

\bibitem[SGA73]{SGA4III73}
\emph{Th\'{e}orie des topos et cohomologie \'{e}tale des sch\'{e}mas. {T}ome
  3}, Lecture Notes in Mathematics, Vol. 305, Springer-Verlag, Berlin-New York,
  1973, S\'{e}minaire de G\'{e}om\'{e}trie Alg\'{e}brique du Bois-Marie
  1963--1964 (SGA 4), Dirig\'{e} par M. Artin, A. Grothendieck et J. L.
  Verdier. Avec la collaboration de P. Deligne et B. Saint-Donat. \MR{0354654}

\bibitem[Vez15]{Vezzosi15}
Gabriele Vezzosi, \emph{A model structure on relative dg-{L}ie algebroids},
  Stacks and categories in geometry, topology, and algebra, Contemp. Math.,
  vol. 643, Amer. Math. Soc., Providence, RI, 2015, pp.~111--118. \MR{3381471}

\bibitem[Vit14]{Vitagliano2014}
Luca Vitagliano, \emph{On the strong homotopy {Lie-Rinehart} algebra of a
  foliation}, Commun. Contemp. Math \textbf{16} (2014), no.~6.

\bibitem[Wei94]{WeibelBook94}
Charles~A. Weibel, \emph{An introduction to homological algebra}, Cambridge
  Studies in Advanced Mathematics, vol.~38, Cambridge University Press,
  Cambridge, 1994. \MR{1269324}

\bibitem[WX91]{WeinsteinX91}
Alan Weinstein and Ping Xu, \emph{Extensions of symplectic groupoids and
  quantization}, J. Reine Angew. Math. \textbf{417} (1991), 159--189.
  \MR{1103911}

\end{thebibliography}
\bibliographystyle{amsalpha}

\end{document}